\documentclass{amsart}

\usepackage{amsmath} 
\usepackage{amssymb}

\newtheorem{theorem}{Theorem}[section] 
\newtheorem{claim}{Claim}[theorem]
\newtheorem{lemma}[theorem]{Lemma} 
\newtheorem{proposition}[theorem]{Proposition} 
\newtheorem{corollary}[theorem]{Corollary} 

\theoremstyle{definition}
\newtheorem{definition}[theorem]{Definition}

\theoremstyle{remark}
\newtheorem{remark}[theorem]{Remark}
\newtheorem{notation}[theorem]{Notation}

\newtheorem{context}[theorem]{Context}

\newcommand{\forces}{\Vdash} 
\newcommand{\bV}{{\bf V}} 
\newcommand{\Cmk}{{\mathbb C}_{\mu,\kappa}}  
\newcommand{\bP}{{\mathbb P}}
\newcommand{\bQ}{{\mathbb Q}}
\newcommand{\dbQ}{{\dot{\mathbb Q}}}

\newcommand{\lesdot}{\mathrel{\mathord{<}\!\!\raise 
0.8 pt\hbox{$\scriptstyle\circ$}}} 

\newcommand{\cov}{{\rm cov}\/}
\newcommand{\gd}{{\mathfrak d}} 
\newcommand{\gdk}{{\mathfrak d}_\kappa}
\newcommand{\gdklm}{{\mathfrak d}_{\kappa,\lambda}^{\kappa,{<}\mu}}
\newcommand{\bgdk}{\bar{\mathfrak d}_\kappa}
\newcommand{\gbk}{{\mathfrak b}_\kappa}

\newcommand{\bairek}{{}^{\textstyle\kappa}\kappa} 
\newcommand{\cantork}{{}^{\textstyle\kappa} 2} 
\newcommand{\seqk}{{}^{\textstyle<\!\kappa}\kappa} 
\newcommand{\rest}{{\restriction}}
\newcommand{\Dom}{{\rm Dom}}
\newcommand{\Rng}{{\rm Rng}}

\newcommand{\Mkk}{{\bf M}_{\kappa,\kappa}}
\newcommand{\Mkkk}{{\bf M}_{\kappa,\kappa}^\kappa}
\newcommand{\cf}{{\rm cf}\/} 
\newcommand{\Gcc}{{\Game^{\rm cc}_{\varepsilon,\kappa}}}
\newcommand{\st}{{\rm st}}
\newcommand{\nst}{\dot{\rm st}}

\newcommand{\cl}{{\rm cl}}
\newcommand{\cA}{{\mathcal A}}
\newcommand{\cC}{{\mathcal C}}
\newcommand{\cH}{{\mathcal H}}
\newcommand{\cF}{{\mathcal F}}
\newcommand{\gF}{{\mathfrak F}}
\newcommand{\cI}{{\mathcal I}}
\newcommand{\cJ}{{\mathcal J}}
\newcommand{\cO}{{\mathcal O}}
\newcommand{\bbP}{{\mathbb P}}
\newcommand{\cP}{{\mathcal P}}
\newcommand{\cX}{{\mathcal X}}
\newcommand{\cY}{{\mathcal Y}}

\newcommand{\Pkl}{{\mathcal P}_\kappa(\lambda)}
\newcommand{\Pm}{{\mathcal P}_\mu}

\newcommand{\lh}{{\rm lh}}
\newcommand{\conc}{{}^\frown\!}

\newcommand{\NSk}{{\mathcal N\! S}_\kappa}
\newcommand{\NSkl}{{\mathcal N\! S}^\kappa_{\kappa,\lambda}}
\newcommand{\NS}{{\mathcal N\! S}}
\newcommand{\cof}{{\rm cof\/}}
\newcommand{\add}{{\rm add\/}}

\newcommand{\bcof}{\overline{\rm cof}\/}

\begin{document}

\title{Cofinality of the nonstationary ideal}

\author{Pierre Matet}
\address{Math\'ematiques\\ 
Universit\'e de Caen\\
BP 5186\\
14032 Caen, Cedex\\ 
France} 
\email{Matet@math.unicaen.fr}

\author{Andrzej Ros{\l}anowski}
\address{Department of Mathematics\\
 University of Nebraska at Omaha\\
 Omaha, NE 68182-0243, USA
 and Department of Mathematics\\
 University of Northern Iowa\\
 Cedar Falls, IA 50614-0506, USA} 
\email{roslanow@member.ams.org}
\urladdr{http://www.unomaha.edu/$\sim$aroslano}

\author{Saharon Shelah}
\address{Institute of Mathematics\\
 The Hebrew University of Jerusalem\\
 91904 Jerusalem, Israel\\
 and  Department of Mathematics\\
 Rutgers University\\
 New Brunswick, NJ 08854, USA}
\email{shelah@math.huji.ac.il}
\urladdr{http://www.math.rutgers.edu/$\sim$shelah}
\thanks{The research of the third author was partially supported by the
 Israel Science Foundation. Publication 799} 


\begin{abstract}
We show that the reduced cofinality of the nonstationary ideal $\NSk$ on a
regular uncountable cardinal $\kappa$ may be less than its cofinality, where
the reduced cofinality of $\NSk$ is the least cardinality of any family
$\cF$ of nonstationary subsets of $\kappa$ such that every nonstationary
subset of $\kappa$ can be covered by less than $\kappa$ many members of
$\cF$. 
\end{abstract}

\maketitle

\setcounter{section}{-1}

\section{Introduction}
Let $\kappa$ be a regular uncountable cardinal. For $C\subseteq\kappa$ and
$\gamma\leq\kappa$, we say that $\gamma$ is {\em a limit point of of $C$} if
$\bigcup(C\cap\gamma)=\gamma>0$. $C$ is {\em closed unbounded\/} if $C$ is a
cofinal subset of $\kappa$ containing all its limit points less than
$\kappa$, A set $A\subseteq\kappa$ is {\em nonstationary\/} if $A$ is
disjoint from some closed unbounded subset $C$ of $\kappa$. The
nonstationary subsets of $\kappa$ form an ideal on $\kappa$ denoted by
$\NSk$. The {\em cofinality\/} of this ideal, $\cof(\NSk)$, is the least
cardinality of any family $\cF$ of nonstationary subsets of $\kappa$ such
that every nonstationary subset of $\kappa$ is contained in a member of
$\cF$. The {\em reduced cofinality\/} of $\NSk$, $\bcof(\NSk)$, is the least
cardinality of any $\cF\subseteq\NSk$ such that every nonstationary subset
of $\kappa$ can be covered by less than $\kappa$ many members of $\cF$. This
paper addresses the question whether $\bcof(\NSk)=\cof(\NSk)$. Note that 
\[\kappa^+\leq\bcof(\NSk)\leq\cof(\NSk)\leq2^\kappa,\]
so under GCH we have $\bcof(\NSk)=\cof(\NSk)$.

Let $\cantork$ be endowed with the $\kappa$--box product topology, 2 itself
considered discrete. We say that a set $W\subseteq \cantork$ is
$\kappa$--meager if there is a sequence $\langle U_\alpha:\alpha<\kappa
\rangle$ of dense open subsets of $\cantork$ such that
$W\cap\bigcap\limits_{\alpha<\kappa} U_\alpha=\emptyset$. The {\em covering
number for the category\/} of the space $\cantork$, denoted $\cov(\Mkk)$, is
the least cardinality of any collection $\cX$ of $\kappa$--meager subsets of
$\cantork$ such that $\bigcup\cX=\cantork$. It can be established that 
\[\cov(\Mkk)\leq\cof(\NSk)\leq\big(\bcof(\NSk)\big)^{<\kappa}.\]
It follows that if $\bcof(\NSk)<\cov(\Mkk)$ and the Singular Cardinals
Hypothesis holds true, then $\cf(\bcof(\NSk))<\kappa$ and $\cof(\NSk)=
(\bcof(\NSk))^+$. We prove:  

\begin{theorem}
\label{thone}
Assume GCH. Then there is a $\kappa$--complete, $\kappa^+$--cc forcing
notion $\bbP$ such that 
\[\forces_{\bbP}\mbox{`` }\bcof(\NSk)=\kappa^{+\omega}\ \mbox{ and }\
\cof(\NSk)=\kappa^{(\omega+1)}\mbox{ ''}.\]
\end{theorem}

\noindent What about the consistency of ``$\bcof(\NSk)$ is regular and 
$\bcof(\NSk)<\cov(\Mkk)$''? We establish: 

\begin{theorem}
\label{thtwo}
It is consistent, relative to the existence of a cardinal $\nu$ such that
$o(\nu)=\nu^{++}$, that $\bcof(\NS_{\omega_1})=\aleph_{\omega+1}$ and
$\cov({\bf M}_{\omega_1,\omega_1})=\aleph_{\omega+2}$.
\end{theorem}

The structure of the paper is as follows. In Section 1, for each infinite
cardinal $\mu\leq\kappa$ we introduce the ${<}\mu$--cofinality
$\bcof^{<\mu}(\NSk)$ and the ${<}\mu$--dominating number
$\gd^{<\mu}_\kappa$, and we show that these two numbers are equal. Section 2
is concerned with a variant of $\gd^{<\mu}_\kappa$ denoted by
$\gd^{\cl,<\mu}_\kappa$ (where $\cl$ stands for ``club''). We establish that
$\gd^{<\mu}_\kappa=\gd^{\cl,<\mu}_\kappa$ if $\mu>\omega$.

$\NSk$ is the smallest normal ideal on $\kappa$. Section 3 deals with
$\NS^\kappa_{\kappa,\lambda}$, the smallest $\kappa$--normal ideal on
$\Pkl$. We compute $\cof^{<\mu}(\NS^\kappa_{\kappa,\lambda})$ and give
examples of situations when $\cof^{<\mu}(\NS^\kappa_{\kappa,\lambda})
<\cof(\NS^\kappa_{\kappa,\lambda})$. 

In the following section we present some basic facts regarding the ideal of
$\kappa$--meager subsets of $\cantork$ and its covering number. 

The final three sections of the paper present the consistency results
mentioned in Theorems \ref{thone}, \ref{thtwo} above. First, in Section 5 we
introduce {\em manageability}, a property of ${<}\kappa$--complete
$\kappa^+$--cc forcing notions which implies preservation of non-meagerness
of subsets of $\bairek$ and which can be iterated. Next, in Section 6, we
define one-step forcing and verify that it has all required properties. The
final section gives the applications obtained by iterating this forcing
notion.

\begin{notation}
\label{notation}
Our notation is rather standard and compatible with that of classical
textbooks (like Jech \cite{J}). In forcing we keep the older (Cohen's)
convention that {\em a stronger condition is the larger one}. Some of our
conventions are listed below.   
\begin{enumerate}
\item For a forcing notion $\bP$, $\Gamma_\bP$ stands for the canonical
$\bP$--name for the generic filter in $\bP$. With this one exception, all
$\bP$--names for objects in the extension via $\bP$ will be denoted with a
dot above (e.g.~$\dot{\tau}$, $\dot{X}$). The weakest element of $\bP$ will
be denoted by $\emptyset_\bP$ (and we will always assume that there is one,
and that there is no other condition equivalent to it). In iterations, if
$\bar{\bQ}=\langle\bP_\zeta,\dot{\bQ}_\zeta:\zeta<\zeta^*\rangle$ and
$p\in\lim(\bar{\bQ})$, then we keep convention that $p(\alpha)=
\dot{\emptyset}_{\dot{\bQ}_\alpha}$ for $\alpha\in\zeta^* \setminus
\Dom(p)$.  
\item Ordinal numbers will be denoted by $\alpha,\beta,\gamma,\delta,
\varepsilon,\zeta,\xi$ and also by $i,j$ (with possible sub- and
superscripts).\\  
Infinite cardinal numbers will be called $\theta,\iota,\mu,\nu,\tau$ (with
possible sub- and superscripts); $\kappa$ is our fixed regular uncountable
cardinal,  $\lambda$ will denote a fixed cardinal $>\kappa$ (in Section 3).  
\item By $\chi$ we will denote a {\em sufficiently large} regular cardinal
and by $\cH(\chi)$ the family of all sets hereditarily of size less than
$\chi$. Moreover, we fix a well ordering $<^*_\chi$ of $\cH(\chi)$. 
\item A bar above a letter denotes that the object considered is a sequence;
usually $\bar{X}$ will be $\langle X_i:i<\zeta\rangle$, where $\zeta$
denotes the length of $\bar{X}$. For a set $A$ and a cardinal $\mu$, the
set of all sequences of members of $A$ of length $\mu$ (length $<\mu$,
respectively), will be denoted by ${}^{\textstyle\mu}A$ (${}^{\textstyle
{<}\mu}A$, respectively). 
\end{enumerate}
\end{notation}

\section{$\cof^{<\mu}(\NSk)$}

\begin{definition}
\label{def1.1}
\begin{enumerate}
\item For a set $A$ and a cardinal $\mu$, $\Pm(A)=\{a\subseteq A:|a|<\mu\}$.
\item Given two infinite cardinals $\mu\leq\tau$, $u(\mu,\tau)$ is the least 
cardinality of any $A\subseteq\Pm(\tau)$ such that $\Pm(\tau)=
\bigcup\limits_{a\in A}\cP(a)$.
\end{enumerate}
\end{definition}

\begin{definition}
\label{def1.2}
Let $S$ be an infinite set and $\cJ$ be an ideal on $S$.
\begin{enumerate}
\item $\cof(\cJ)$ is the least cardinality of any $\cX\subseteq\cJ$ such
that for every $A\in\cJ$, there is $B\in\cX$ with $A\subseteq B$.
\item $\add(\cJ)$ is the least cardinality of any $\cX\subseteq\cJ$ such
that $\bigcup\cX\notin \cJ$.
\item For an infinite cardinal $\mu\leq\add(\cJ)$, $\cof^{<\mu}(\cJ)$ is the
least cardinality of a family $\cX\subseteq\cJ$ such that for every
$A\in\cJ$, there is $\cY\in\Pm(\cX)$ such that $A\subseteq\bigcup\cY$.
\item We let $\bcof(\cJ)=\cof^{<\add(\cJ)}(\cJ)$.
\end{enumerate}
\end{definition}

The following proposition collects some trivialities.

\begin{proposition}
\label{prop1.3}
Let $S$ be an infinite set and $\cJ$ be an ideal on $S$. Then:
\begin{enumerate}
\item[(i)]    $\cof^{<\omega}(\cJ)=\cof(\cJ)$.
\item[(ii)]   If $\mu,\nu$ are two infinite cardinals with $\mu\leq
\nu\leq\add(\cJ)$,\\
then $\cof^{<\nu}(\cJ)\leq\cof^{<\mu}(\cJ)$.
\item[(iii)]  $\cof(\cJ)\leq u(\mu,\cof^{<\mu}(\cJ))$ for every infinite
cardinal $\mu\leq\add(\cJ)$. 
\item[(iv)]   $\add(\cJ)\leq\bcof(\cJ)$.
\end{enumerate}
\end{proposition}

The following is well-known (see, e.g., Matet, P\'ean and Shelah
\cite{MPSh:713}): 

\begin{lemma}
\label{lem1.5}
Let $\mu$ be a regular infinite cardinal. Then $u(\mu,\mu^{+n})=\mu^{+n}$
for every $n<\omega$.
\end{lemma}

\begin{proposition}
\label{prop1.6}
Let $S$ be an infinite set and $\cJ$ be an ideal on $S$ such that
$\big(\add(\cJ)\big)^{+\omega}\leq\cof(\cJ)$. Then $\big(\add(\cJ)\big
)^{+\omega}\leq\bcof(\cJ)$. 
\end{proposition}

\begin{proof}
Use Lemma \ref{lem1.5}.
\end{proof}

With these preliminaries out of the way, we can concentrate on ideals on
$\kappa$. If there is a family of size $\kappa^{+\omega}$ of pairwise almost
disjoint cofinal subsets of $\kappa$, then there is a $\kappa$--complete
ideal $\cJ$ on $\kappa$ such that $\bcof(\cJ)<\cof(\cJ)$ (see Matet and
Pawlikowski \cite{MPxx}). Also the converse holds. 

\begin{proposition}
\label{prop1.7}
Let $\cJ$ be a $\kappa$--complete ideal on $\kappa$ such that $\bcof(\cJ)>
\kappa$. Then there is a family $Q\subseteq \cJ\setminus\cP_\kappa(\kappa)$
such that $|Q|=\bcof(\cJ)$ and $A\cap B\in\cP_\kappa (\kappa)$ for any two
distinct members $A$ and $B$ of $Q$. 
\end{proposition}

\begin{proposition}
\label{prop1.8}
Suppose $\cJ$ is a normal ideal on $\kappa$ and $\kappa$ is a limit
(regular) cardinal. Then $\bcof(\cJ)=\cof^{<\mu}(\cJ)$ for some infinite
cardinal $\mu<\kappa$.
\end{proposition}

\begin{proof}
Assume that the conclusion fails. Fix $\cX\subseteq\cJ$ such that $|\cX|=
\bcof(\cJ)$ and 
\[\cJ=\bigcup\{\cP({\textstyle \bigcup} X):X\in\cP_\kappa(\cX)\}.\]
Set $\cY=\{A\cup\beta: A\in \cX\ \&\ \beta\in\kappa\}$. Note that $|\cY|=
\bcof(\cJ)$. Assume that for each infinite cardinal $\mu<\kappa$ we may
select a set $B_\mu\in\cJ$ so that $B_\mu\nsubseteq\bigcup Y$ for any
$Y\in\cP_\mu(\cY)$. Now let $B$ be the set of all $\alpha<\kappa$ such that
$\alpha\in B_\mu$ for some infinite cardinal $\mu<\alpha$. Since $B\in\cJ$
(by normality of $\cJ$), there must be $X\in\cP_\kappa(\cX)$ such that $B
\subseteq\bigcup X$. Let $\tau$ be any infinite cardinal such that
$|X|<\tau<\kappa$. Then $B_\tau\subseteq\bigcup\limits_{A\in X}(A\cup
(\tau+1))$, which is a contradiction. 
\end{proof}

Arguing as in Proposition \ref{prop1.8}, we get:

\begin{proposition}
\label{prop1.9}
Suppose $\cJ$ is a $\kappa$--complete ideal on $\kappa$ and $\nu$ is an
uncountable limit cardinal $<\kappa$. Then there is an infinite cardinal
$\mu<\nu$ such that $\cof^{<\nu}(\cJ)=\cof^{<\mu}(\cJ)$. Moreover, the least
such $\mu$ is either $\omega$, or a successor cardinal. 
\end{proposition}

The remainder of this section is concerned with $\cof^{<\mu}(\NSk)$. Let us
recall the definition of the bounding number $\gbk$:

\begin{definition}
\label{def1.10}
The {\em $\kappa$--bounding number\/} $\gbk$ is the least cardinality of any
$\cF\subseteq\bairek$ with the property that for every $g\in\bairek$, there
is $f\in\cF$ such that 
\[|\{\alpha<\kappa:g(\alpha)\leq f(\alpha)\}|=\kappa.\]
\end{definition}

The following is proved in Matet and Pawlikowski \cite{MPxx}:

\begin{proposition}
\label{prop1.11}
\begin{enumerate}
\item[(i)]  $\bcof(\NSk)\geq\gbk$.
\item[(ii)] If $\bcof(\NSk)=\gbk$, then $\bcof(\NSk)=\cof(\NSk)$.
\end{enumerate}
\end{proposition}

\begin{proposition}
\label{prop1.12}
Let $\mu$ be an infinite cardinal $\leq\kappa$. Then 
\[\mbox{either }\ \cf(\cof^{<\mu}(\NSk))<\mu,\ \mbox{ or }\ \cf(
\cof^{<\mu}(\NSk))\geq\gbk.\]  
\end{proposition}

\begin{proof}
Suppose to the contrary that $\mu\leq\cf(\cof^{<\mu}(\NSk))=\tau<
\gbk$. For $\alpha<\tau$ select $\cX_\alpha\subseteq\NSk$ so that  
\begin{enumerate}
\item[(i)]   $|\cX_\alpha|<\cof^{<\mu}(\NSk)$,
\item[(ii)]  $\cX_\beta\subseteq \cX_\alpha$ for $\beta<\alpha$, 
\item[(iii)] $\NSk=\bigcup\{\cP(\bigcup X):X\in\cP_\mu(
\bigcup\limits_{\alpha<\tau}\cX_\alpha)\}$. 
\end{enumerate}
For $\alpha<\tau$, set $\cY_\alpha=\{A\cup\beta:A\in\cX_\alpha\ \&\
\beta\in\kappa\}$ and pick $B_\alpha\in\NSk$ so that $B_\alpha
\nsubseteq\bigcup Y$ for any $Y\in \cP_\mu(\cY_\alpha)$. By a result of
Balcar and Simon (see \cite[Theorem 5.25]{BS89}), there is $B\in\NSk$ such
that $|B_\alpha\setminus B|<\kappa$ for every $\alpha<\tau$. Select $X\in
\cP_\mu(\bigcup\limits_{\alpha<\tau}\cX_\alpha)$ so that $B\subseteq \bigcup
X$. There is $\gamma<\tau$ such that $X\subseteq\cX_\gamma$. Then $B_\gamma
\subseteq\bigcup\limits_{A\in X}(A\cup\beta)$ for some $\beta\in\kappa$,
which is a contradiction. 
\end{proof}

\begin{definition}
\label{def1.13}
Let $\tau\leq\kappa$. A family $\cF\subseteq\bairek$ is called
\begin{itemize}
\item {\em a dominating family\/} if 
\[(\forall h\in\bairek)(\exists f\in\cF)(\forall j<\kappa)(h(j)<f(j)),\]
\item {\em a ${<}\tau$--dominating family\/} if 
\[(\forall h\in\bairek)(\exists F\in \cP_\tau(\cF))(\forall j<\kappa)(h(j)<
\sup\{f(j):f\in F\}).\] 
\end{itemize}
We define dominating numbers $\gd_\kappa,\gd_\kappa^{<\tau}$ by:
\[\begin{array}{l}
\gd_\kappa=\min\{|\cF|:\cF\subseteq\bairek\mbox{ is a dominating family
}\},\\  
\gd_\kappa^{<\tau}=\min\{|\cF|:\cF\subseteq\bairek\mbox{ is a
${<}\tau$--dominating family }\}. 
  \end{array}\]
We let $\bgdk=\gdk^{<\kappa}$ and for an infinite cardinal $\mu<\kappa$ we
put $\gdk^\mu=\gdk^{<\mu^+}$.
\end{definition}

Note that $\gdk^{<\omega}=\gdk$. Landver \cite{La90} established that
$\cof(\NSk)=\gdk$. His result can be generalized as follows:

\begin{theorem}
\label{thm1.14}
Let $\mu$ be an infinite cardinal $\leq\kappa$. Then $\cof^{<\mu}(\NSk)=
\gdk^{<\mu}$. 
\end{theorem}

\begin{proof}
Set $\tau=\cof^{<\mu}(\NSk)$. First we establish that $\gdk^{<\mu}\leq
\tau$. Select s family $\cC$ of size $\tau$ of closed unbounded subsets of
$\kappa$ so that for every closed unbounded subset $D$ of $\kappa$, there
is $X\in\cP_\mu(\cC)\setminus\{\emptyset\}$ with $\bigcap X\subseteq D$. For
$U\in\cP_\omega(\cC)\setminus\{\emptyset\}$ define $f_U\in\bairek$ by
$f_U(\alpha)=\min\big(\bigcap U\setminus (\alpha+1)\big)$. Note that
$f_V(\alpha)\leq f_U(\alpha)$ whenever $V\in \cP(U)\setminus\{\emptyset\}$. 
Now given $g\in\bairek$, let $D$ be the set of all limit ordinals $\delta$
such that $0<\delta<\kappa$ and $g(\alpha)<\delta$ for every
$\alpha<\delta$. Pick $X\in\cP_\mu(\cC)\setminus\{\emptyset\}$ so that
$\bigcap X\subseteq D$. Define $h\in\bairek$ by
\[h(\alpha)=\sup\big\{f_U(\alpha):U\in\cP_\omega(X)\setminus\{\emptyset\}
\big\}.\] 
Let $\alpha<\kappa$ and $C\in X$. First, suppose that there is $W\in
\cP_\omega(X)\setminus\{\emptyset\}$ such that $h(\alpha)=f_W(\alpha)$. Then
$h(\alpha)=f_{W\cup\{C\}}(\alpha)$ and hence $h(\alpha)\in C$. Next suppose
that $f_U(\alpha)<h(\alpha)$ for all $U\in\cP_\omega(X)\setminus\{\emptyset
\}$. Then $h(\alpha)$ is a limit ordinal. Set $\iota=\cf(h(\alpha))$ and
pick an increasing sequence $\langle\gamma_\beta:\beta<\iota\rangle$ cofinal
in $h(\alpha)$. For $\beta<\iota$, select $T_\beta\in\cP_\omega(X)\setminus
\{\emptyset\}$ with $\gamma_\beta<f_{T_\beta}(\alpha)$, and set
$\delta_\beta=f_{T_\beta\cup\{C\}}(\alpha)$. Note that $\delta_\beta\in
C$. Obviously, the sequence $\langle\delta_\beta:\beta<\iota\rangle$ is
cofinal in $h(\alpha)$, and consequently $h(\alpha)\in C$. Thus for each
$\alpha<\kappa$, $h(\alpha)$ belongs to $\bigcap X$ and therefore to
$D$. Since clearly $h(\alpha)>\alpha$, it follows that $h(\alpha)>g(\alpha
)$. 

It remains to show that $\gdk^{<\mu}\geq\tau$. Let $\gF$ be the set of all
strictly increasing functions from $\kappa$ to $\kappa$. Select $\cF
\subseteq\gF$ so that 
\begin{enumerate}
\item[(a)] $|\cF|=\gdk^{<\mu}$, and 
\item[(b)] given $g\in\bairek$, there is $F_g\in\cP_\mu(\cF)$ such that 
\[(\forall\alpha<\kappa)(g(\alpha)<\sup\{f(\alpha):f\in F_g\}).\]
\end{enumerate}
For $f\in\gF$, let $C_f$ be the set of all limit ordinals $\alpha$ such that
$0<\alpha<\kappa$ and $f(\beta)<\alpha$ for every $\beta<\alpha$. Easily
\[\NSk=\{A\subseteq\kappa: (\exists g\in\gF)(A\cap C_g=\emptyset)\}\]
(see, e.g., \cite{MPxx}) and (as $\bigcap\limits_{f\in F_g} C_f\subseteq
C_g$ for every $g\in\gF$) it follows that $\tau\leq |\cF|$. 
\end{proof}

By Propositions \ref{prop1.8} and \ref{prop1.9} and Theorem \ref{thm1.14},
to determine the value of $\cof^{<\mu}(\NSk)$ for every infinite cardinal
$\mu\leq\kappa$, it suffices to compute $\gdk$ and $\gdk^\tau$ for every
infinite cardinal $\tau<\kappa$.

\section{$\gdk^{\cl,{<}\mu}$}
It is straightforward to check that $\gdk^{<\mu}$ is the least cardinality
of a family $\cF\subseteq\bairek$ such that 
\[(\forall g\in\bairek)(\exists F\in\cP_\mu(\cF))(|\{\alpha\in\kappa:
g(\alpha)\geq \sup\{f(\alpha):f\in F\}|<\kappa).\]
In this section we discuss the variant that arises if we replace the
noncofinal ideal on $\kappa$ by the nonstationary ideal on $\kappa$. 

\begin{definition}
\label{def2.1}
\begin{enumerate}
\item $\gdk^{\cl}$ is the least cardinality of a family $\cF\subseteq
\bairek$ with the property that for every $g\in\bairek$, there is $f\in\cF$
such that 
\[\{\alpha\in\kappa:g(\alpha)\geq f(\alpha)\}\in\NSk.\]
\item For an infinite cardinal $\mu\leq\kappa$, $\gdk^{\cl,{<}\mu}$ is the
least cardinality of a family $\cF\subseteq\bairek$ with the property that
for every $g\in\bairek$, there is $F\in\cP_\mu(\cF)$ such that 
\[\big\{\alpha\in\kappa:g(\alpha)\geq\sup\{f(\alpha):f\in F\}\big\}\in\NSk.\]
\end{enumerate}
\end{definition}

Note that $\gdk^{\cl,{<}\omega}=\gdk^{\cl}$. It is simple to check that
$\cf(\gdk^{\cl})\geq\gbk$. 

\begin{theorem}
\label{thm2.2}
For every uncountable cardinal $\mu\leq\kappa$, 
\[\gdk^{\cl,{<}\mu}=\gdk^{<\mu}.\]
\end{theorem}

Theorem \ref{thm2.2} easily follows from the next two lemmas. 

\begin{lemma}
\label{lem2.3}
Let $\mu$ be an uncountable limit cardinal $\leq\kappa$. Then
$\gdk^{\cl,{<}\mu}\geq\gdk^{\cl,{<}\tau}$ for some infinite cardinal
$\tau<\mu$. 
\end{lemma}

\begin{proof}
The proof is similar to that of Proposition \ref{prop1.9}. Suppose that the
conclusion fails. Fix a family $\cF\subseteq\bairek$ such that $|\cF|=
\gdk^{\cl,{<}\mu}$ and 
\[(\forall g\in\bairek)(\exists F\in\cP_\mu(\cF))(\big\{\alpha\in\kappa:
g(\alpha)\geq \sup\{f(\alpha):f\in F\}\big\}\in\NSk).\]
For each infinite cardinal $\tau<\mu$ select $g_\tau \in\bairek$ so that for
every $F\in\cP_\tau(\cF)$ we have
\[\big\{\alpha\in\kappa:g_\tau(\alpha)\geq\sup\{f(\alpha):f\in F\}\big\}
\notin\NSk.\] 
Define $g\in\bairek$ so that $g(\alpha)\geq g_\tau(\alpha)$ for every
infinite cardinal $\tau<\mu$ such that $\tau<\alpha$. Now pick
$F\in\cP_\mu(\cF)$ such that 
\[\big\{\alpha\in\kappa:g(\alpha)\geq \sup\{f(\alpha):f\in F\}\big\}\in
\NSk.\] 
Let $\tau$ be any infinite cardinal with $|F|<\tau<\mu$. Obviously, $F\in
\cP_\tau(\cF)$ and 
\[\big\{\alpha\in\kappa:g_\tau(\alpha)\geq \sup\{f(\alpha):f\in F\}\big\}\in
\NSk,\]
a contradiction.
\end{proof}

To establish the following lemma, we adapt the proof of Theorem 5 in
Cummings and Shelah \cite{CuSh:541}.

\begin{lemma}
\label{lem2.4}
Let $\mu$ be a regular uncountable cardinal $\leq\kappa$. Then $\gdk^{<\mu}
\leq\gdk^{\cl,{<}\mu}$. 
\end{lemma}

\begin{proof}
Select a family $\cF\subseteq\bairek$ such that
\begin{enumerate}
\item[(a)] every member of $\cF$ is increasing, 
\item[(b)] $|\cF|=\gdk^{\cl,{<}\mu}$, and
\item[(c)] for each $g\in\bairek$, there is $F\in\cP_\mu(\cF)$ such that
\[\big\{\alpha\in\kappa:g(\alpha)\geq\sup\{f(\alpha):f\in F\}\big\}\in
\NSk.\]
\end{enumerate}
Now fix $g\in\bairek$. Stipulate that $g_{-1}=g$. By induction on
$n\in\omega$ choose a closed unbounded subset $C_n$ of $\kappa$,
$g_n\in\bairek$, $h_n\in\bairek$ and $F_n\in\cP_\mu(\cF)$ so that 
\begin{enumerate}
\item[(i)]    $C_{n+1}\subseteq C_n$,
\item[(ii)]   $g_{n-1}(\alpha)<\sup\{f(\alpha):f\in F_n\}$ for all
$\alpha\in C_n$,
\item[(iii)]  $h_n(\beta)=\min(C_n\setminus (\beta+1))$, 
\item[(iv)]   $g_n(\beta)=\sup\big(\Rng(g_{n-1}\rest (h_n(\beta)+1))\big)$.  
\end{enumerate}
Note that, by (iii) and (iv), $g(\beta)\leq g_0(\beta)\leq g_1(\beta)\leq
\ldots$ for all $\beta\in\kappa$. Set $F=\bigcup\limits_{n\in\omega} F_n$
and $\zeta=\sup\{\min(C_n):n\in\omega\}$. We are going to show that
$g(\gamma)<\sup\{f(\gamma): f\in F\}$ whenever $\zeta<\gamma<\kappa$. To
this end suppose that $\zeta<\gamma<\kappa$. By (i), there are $m\in\omega$
and $\xi\in\kappa$ such that $\xi=\sup(\gamma\cap C_n)$ whenever $m\leq
n<\omega$. By (iii), $h_m(\xi)\geq\gamma$ and so (by (iv)) $g(\gamma)\leq
g_{m-1}(\gamma)\leq g_m(\xi)$. Since $\gamma>\zeta$ we also have $\gamma\cap
C_{m+1}\neq\emptyset$. Hence $\xi\in C_{m+1}$ and consequently, by (ii), 
\[g(\gamma)\leq g_m(\xi)<\sup\{f(\xi):f\in F_{m+1}\}\leq\sup\{f(\gamma):f\in
F_{m+1}\}.\] 
\end{proof}
 
Theorem \ref{thm2.2} implies that $\gdk^\omega\leq\gdk^{\cl}\leq\gdk$. We
mention that it was shown in Cummings and Shelah \cite{CuSh:541} that
$\gdk^{\cl}=\gdk$ if $\kappa>\beth_\omega$.

\section{$\cof^{{<}\mu}(\NS^\kappa_{\kappa,\lambda})$}
Throughout this section $\lambda$ denotes a fixed cardinal $>\kappa$. Our
object of study will be the ideal $\NSkl$, a $\cP_\kappa(\lambda)$ version
of $\NSk$. 

\begin{definition}
\label{def3.1}
For a regular uncountable cardinal $\nu$ and a cardinal $\tau\geq\nu$,
$\cJ_{\nu,\tau}$ is the set of all $A\subseteq\cP_\nu(\tau)$ such that for
some $a\in\cP_\nu(\tau)$ we have $\{b\in A:a\subseteq b\}=\emptyset$. 
\end{definition}

It is straightforward to check that $\cJ_{\nu,\tau}$ is a $\nu$--complete
ideal on $\cP_\nu(\tau)$. 

\begin{definition}
\label{def3.2}
\begin{enumerate}
\item An ideal $\cJ$ of $\cP_\kappa(\lambda)$ is $\kappa$--normal if given
$A\in\cJ^+$ and $f:A\longrightarrow\kappa$ such that $f(a)\in a\cap\kappa$ for
all $a\in A$, there is $B\in\cJ^+\cap\cP(A)$ such that $f$ is constant on
$B$. 
\item The smallest $\kappa$--normal ideal $\cP_\kappa(\lambda)$ containing
$\cJ_{\kappa,\lambda}$ is denoted by $\NSkl$.
\item For $f\in{}^\kappa(\cP_\kappa(\lambda))$ we let
\[C_f\stackrel{\rm def}{=}\{a\in\cP_\kappa(\lambda):a\cap\kappa\neq
\emptyset\ \mbox{ and }\ \bigcup_{\alpha\in a\cap\kappa} f(\alpha)\subseteq
a\}.\] 
\end{enumerate}
\end{definition}

The following lemma is due to Abe.

\begin{lemma}
[Abe {\cite{Ab97}}]
\label{lem3.4}
Let $A\subseteq\cP_\kappa(\lambda)$. Then 
\[A\in\NSkl\quad\mbox{ if and only if }\quad (\exists f\in {}^\kappa(
\cP_\kappa(\lambda)))(A\cap C_f=\emptyset).\]
\end{lemma}

Our purpose in this section is to compute the value of
$\cof^{<\mu}(\NSkl)$. We will need an analogue of $\gdk^{<\mu}$ defined in
\ref{def3.5}(1) below. 

\begin{definition}
\label{def3.5}
Let $\mu\leq\kappa$ be an infinite cardinal.
\begin{enumerate}
\item $\gdklm$ is the least cardinality of a family $\cX$ of functions from
$\kappa$ to $\cP_\kappa(\lambda)$ with the property that
\[\big(\forall g\in {}^\kappa(\cP_\kappa(\lambda))\big)\big(\exists X\in
\cP_\mu(\cX)\big)\big(\forall\alpha\in\kappa\big)\big(g(\alpha)\subseteq
\bigcup_{f\in X}f(\alpha)\big).\]
\item $\cov(\lambda,\kappa^+,\kappa^+,\mu)$ is the least cardinality of a
family $\cX\subseteq\cP_{\kappa^+}(\lambda)$ such that 
\[\big(\forall B\in\cP_{\kappa^+}(\lambda)\big)\big(\exists X\in
\cP_\mu(\cX)\big)\big(B\subseteq\bigcup X\big).\]
\end{enumerate}
\end{definition}

\begin{theorem}
\label{thm3.7}
Let $\mu$ be an infinite cardinal $\leq\kappa$. Then 
\[\cof^{<\mu}(\NSkl)=\gdklm=\max\{\gdk^{<\mu},\cov(\lambda,\kappa^+,
\kappa^+,\mu)\}.\]
\end{theorem}

Theorem \ref{thm3.7} is an immediate consequence of Lemmas
\ref{lem3.8}--\ref{lem3.11} below.

\begin{lemma}
\label{lem3.8}
Let $\mu$ be an infinite cardinal $\leq\kappa$. Then 
\[\cov(\lambda,\kappa^+,\kappa^+,\mu)\leq \cof^{<\mu}(\NSkl).\]
\end{lemma}

\begin{proof}
By \ref{lem3.4} we may pick a family $\cX\subseteq{}^\kappa(\cP_\kappa(
\lambda))$ with the property that $|\cX|=\cof^{<\mu}(\NSkl)$ and for every
function $g:\kappa \longrightarrow\cP_\kappa(\lambda)$ there is
$X\in\cP_\mu(\cX)$ such that $\bigcap\limits_{f\in X} C_f\subseteq C_g$. For
$f\in \cX$, let $B_f=\kappa\cup\bigcup\limits_{\alpha<\kappa}f(\alpha)\in
\cP_{\kappa^+}(\lambda)$.  

Suppose now that $B\in\cP_{\kappa^+}(\lambda)$. Pick a function $g:\kappa
\longrightarrow\cP_\kappa(\lambda)$ such that $B\subseteq\bigcup\limits_{
\alpha<\kappa}g(\alpha)$. There is $X\in\cP_\mu(\cX)$ such that
$\bigcap\limits_{f\in X} C_f\subseteq C_g$. We are going to show that
$B\subseteq \bigcup\limits_{f\in X}B_f$. To this end suppose
$\alpha<\kappa$ and let us argue that $g(\alpha)\subseteq
\bigcup\limits_{f\in X}B_f$. For $n<\omega$ let $a_n\in\cP_\kappa(
\bigcup\limits_{f\in X} B_f)$ be defined by 
\[a_0=\{\alpha\},\quad\mbox{ and }\quad a_{n+1}=a_n\cup\bigcup_{f\in X}
\bigcup_{\beta\in a_n\cap\kappa} f(\beta),\]
and let $a=\bigcup\limits_{n<\omega} a_n$. Then $\alpha\in a\in
\bigcap\limits_{f\in X} C_f\subseteq C_g$ and consequently
$g(\alpha)\subseteq a\subseteq\bigcup\limits_{f\in X} B_f$.
\end{proof}

\begin{lemma}
\label{lem3.9}
Let $\mu$ be an infinite cardinal $\leq\kappa$. Then 
$\gdk^{<\mu}\leq\cof^{<\mu}(\NSkl)$.
\end{lemma}

\begin{proof}
By Theorem \ref{thm1.14}, it suffices to establish that
$\cof^{<\mu}(\NSk)\leq\cof^{<\mu}(\NSkl)$. Let a family $\cX\subseteq
{}^\kappa(\cP_\kappa(\lambda))$ be such that $|\cX|=\cof^{<\mu}(\NSkl)$
and 
\[\big(\forall B\in\NSkl\big)\big(\exists X\in\cP_\mu(\cX)\setminus\{
\emptyset\}\big)\big(B\cap\bigcap_{f\in X}C_f=\emptyset\big).\]
For $f\in\cX$, let $Z_f$ be the set of all limit ordinals $\alpha>0$ such
that 
\[(\forall\beta<\alpha)(f(\beta)\cap\kappa\subseteq\alpha).\]
Plainly, $Z_f$ is a closed unbounded subset of $\kappa$.  Now given a closed
unbounded subset $T$ of $\kappa$, set $B_T=\{a\in\cP_\kappa(\lambda):a\cap
\kappa\notin T\}$. A simple argument (see, e.g., \cite{MPSh:713}) shows that
$B_T\in\NSkl$. Hence there is $X_T\in\cP_\mu(\cX)\setminus\{\emptyset\}$  
such that $B_T\cap\bigcap\limits_{f\in X_T} C_f=\emptyset$. We will show
that $\bigcap\limits_{f\in X_T} Z_f\subseteq T$. Thus let
$\alpha\in\bigcap\limits_{f\in X_T} Z_f$. Setting $a=\alpha\cup 
\bigcup\limits_{f\in X_T}\bigcup\limits_{\beta<\alpha} f(\beta)$, it is
easy to see that $a\cap\kappa=\alpha$, and so $a\in\bigcap\limits_{f\in X_T}
C_f$. It follows that $\alpha=a\cap\kappa\in T$.
\end{proof}

\begin{lemma}
\label{lem3.10}
Let $\mu$ be an infinite cardinal $\leq\kappa$. Then $\cof^{<\mu}(\NSkl)\leq
\gdklm$. 
\end{lemma}

\begin{proof}
The inequality easily follows from the following observation. 

Suppose $h:\kappa \longrightarrow\cP_\kappa(\lambda)$ and $X\in \cP_\mu\big( 
{}^\kappa(\cP_\kappa(\lambda))\big)$ are such that 
\[\big(\forall\alpha<\kappa\big)\big(h(\alpha)\subseteq\bigcup\limits_{f\in
X} f(\alpha)\big).\]
Then $\bigcap\limits_{f\in X}C_f\subseteq C_h$.  
\end{proof}

\begin{lemma}
\label{lem3.11}
Let $\mu$ be an infinite cardinal $\leq\kappa$. Then 
\[\gdklm\leq\max\{\gdk^{<\mu},\cov(\lambda,\kappa^+,\kappa^+,\mu)\}.\]  
\end{lemma}

\begin{proof}
Fix $\cF$ so that $|\cF|=\gdk^{<\mu}$ and 
\[\big(\forall h\in\bairek\big)\big(\exists F\in\cP_\mu(\cF)\big)\big(
\forall\alpha<\kappa\big)\big(h(\alpha)\leq\sup\{f(\alpha):f\in F\}\big).\] 
Also, fix $\cX\subseteq\cP_{\kappa^+}(\lambda)$ such that $|\cX|=\cov(
\lambda,\kappa^+,\kappa^+,\mu)$ and 
\[\cP_{\kappa^+}(\lambda)=\bigcup\{\cP({\textstyle\bigcup} X):X\in\cP_\mu(
\cX)\}.\] 
For each $a\in\cX$, select a mapping $\phi_a:\kappa\stackrel{\rm onto}{
\longrightarrow}a$. Now for $f\in\cF$ and $a\in\cX$ define $g_{f,a}:\kappa
\longrightarrow\cP_\kappa(\lambda)$ by 
\[g_{f,a}(\alpha)=\{\phi_a(\delta):\delta<f(\alpha)\}.\]
Suppose now that $g:\kappa\longrightarrow\cP_\kappa(\lambda)$. By the choice
of $\cX$, there is $X\in\cP_\mu(\cX)$ such that $\bigcup\limits_{\alpha<
\kappa} g(\alpha)\subseteq\bigcup X$. Choose $h\in\bairek$ such that 
\[\big(\forall\alpha<\kappa\big)\big(g(\alpha)\subseteq\bigcup_{a\in X}\{
\phi_a(\xi):\xi<h(\alpha)\}\big).\] 
Next pick $F\in\cP_\mu(\cF)$ such that $(\forall\alpha<\kappa)(h(\alpha)\leq
\sup\{f(\alpha):f\in F\})$. Then 
\[\big(\forall\alpha<\kappa\big)\big(g(\alpha)\subseteq \bigcup_{f\in F}
\bigcup_{a\in X} g_{f,a}\big).\] 
\end{proof}

Another formula worth noting is:
\[\cof^{<\mu}(\NSkl)=\max\{\cof^{<\mu}(\NSk),\cof^{<\mu}(\cJ_{\kappa^+,
\lambda})\}.\]
This identity follows from Theorems \ref{thm1.14} and \ref{thm3.7} and the
next proposition.

\begin{proposition}
\label{prop3.12}
Let $\mu$ be an infinite cardinal $\leq\kappa$. Then $\cov(\lambda,\kappa^+,
\kappa^+,\mu)=\cof^{<\mu}(\cJ_{\kappa^+,\lambda})$.
\end{proposition}

\begin{proof}
The result easily follows from the following observation. 

Suppose that $\cX\subseteq\cP_{\kappa^+}(\lambda)$ and
$X\in\cP_\mu(\cX)$. Then  
\[\bigcap_{a\in X}\{c\in\cP_{\kappa^+}(\lambda): a\subseteq c\}=\{c\in\cP_{
\kappa^+}(\lambda): {\textstyle\bigcup} X\subseteq c\},\]
and therefore for each $b\in\cP_{\kappa^+}(\lambda)$
\[b\subseteq{\textstyle\bigcup} X\quad\mbox{ if and only if }\quad
\bigcap_{a\in X}\{c\in\cP_{\kappa^+}(\lambda):a\subseteq c\}\subseteq\{c\in  
\cP_{\kappa^+}(\lambda): b\subseteq c\}.\]  
\end{proof}

We next consider special cases when $\cof^{<\mu}(\NSkl)<\cof(\NSkl)$. 

\begin{lemma}
\label{lem3.13}
Let $\mu$ be an infinite cardinal $\leq\kappa$. Then $\cov(\lambda,\kappa^+,
\kappa^+,\mu)\geq\lambda$. 
\end{lemma}

\begin{proof}
It is shown in Matet, P\'ean and Shelah \cite{MPSh:713} that
$\bcof(\cJ_{\kappa^+,\lambda})\geq \lambda$. Now observe that (by
Proposition \ref{prop3.12}) $\cov(\lambda,\kappa^+,\kappa^+,\mu)\geq
\bcof(\cJ_{\kappa^+,\lambda})$. 
\end{proof}

\begin{lemma}
\label{lem3.14}
Suppose $\lambda$ is singular and $\mu$ is a cardinal such that
$\cf(\lambda)<\mu\leq\kappa$. Then $\cov(\lambda,\kappa^+,\kappa^+,\mu)\leq
\sup\{u(\kappa^+,\nu):\kappa<\nu<\lambda\}$. 
\end{lemma}

\begin{proof}
Let $\langle\lambda_\xi:\xi<\cf(\lambda)\rangle$ be an increasing sequence
of cardinals cofinal in $\lambda$. Then, for every
$a\in\cP_{\kappa^+}(\lambda)$, $a=\bigcup\limits_{\xi<\cf(\lambda)}a\cap
\lambda_\xi$. The desired inequality follows.
\end{proof}

\begin{proposition}
\label{prop3.15}
Let $\mu$ be an uncountable cardinal $\leq\kappa$. Then 
\[\cof^{<\mu}\big(\NS^\kappa_{\kappa,\kappa^{+\omega}}\big)=\max\{
\gdk^{<\mu},\kappa^{+\omega}\}.\]
\end{proposition}

\begin{proof}
By Lemmas \ref{lem1.5}, \ref{lem3.13} and \ref{lem3.14} we have
\[\cov(\kappa^{+\omega},\kappa^+,\kappa^+,\mu)=\kappa^{+\omega},\]
so the result follows from Theorem \ref{thm3.7}.
\end{proof}

Thus, if $\gdk^{<\omega_1}\leq\kappa^{+\omega}$, then 
\[\cof^{<\omega_1}(\NS_{\kappa,\kappa^{+\omega}}^\kappa)<\cof(\NS_{
\kappa,\kappa^{+\omega}}^\kappa).\]

\begin{lemma}
\label{lem3.16}
Assume the Singular Cardinals Hypothesis. If $\lambda\geq 2^\kappa$, then 
\[u(\kappa^+,\lambda)=\left\{\begin{array}{ll}
\lambda^+&\mbox{ if }\cf(\lambda)\leq\kappa,\\
\lambda  &\mbox{ otherwise.}
	               	     \end{array}\right.\]
\end{lemma}

\begin{proof}
By Lemma \ref{lem3.13}, $\lambda\leq u(\kappa^+,\lambda)\leq
\lambda^\kappa$. It follows immediately that $u(\kappa^+,\lambda)=\lambda$
if $\cf(\lambda)\geq\kappa^+$. For the other case, use the well--known fact
(see, e.g., \cite{MPSh:713}) that $\cf(u(\kappa^+,\lambda))\geq\kappa^+$.  
\end{proof}

\begin{proposition}
\label{prop3.17}
Assume the Singular Cardinals Hypothesis. If $\lambda\geq 2^\kappa$ and
$\aleph_0\leq\mu\leq\kappa$, then 
\[\cof^{<\mu}(\NSkl)=\left\{\begin{array}{ll}
\lambda^+&\mbox{ if }\mu\leq\cf(\lambda)\leq\kappa,\\
\lambda  &\mbox{ otherwise.}
	               	     \end{array}\right.\]
\end{proposition}

\begin{proof}
By Lemma \ref{lem3.13}, $\cov(\lambda,\kappa^+,\kappa^+,\mu)\geq\lambda
\geq\gdk\geq\gdk^{<\mu}$, so by Theorem \ref{thm3.7}
\[\cof^{<\mu}(\NSkl)=\cov(\lambda,\kappa^+,\kappa^+,\mu).\]

\noindent{\sc Case:}\quad $\cf(\lambda)>\kappa$.\\
By Lemmas \ref{lem3.13} and \ref{lem3.16} we have 
\[\lambda\leq\cov(\lambda,\kappa^+,\kappa^+,\mu)\leq u(\kappa^+,\lambda)
\leq\lambda,\]
and hence $\cov(\lambda,\kappa^+,\kappa^+,\mu)=\lambda$.
\medskip

\noindent{\sc Case:}\quad $\mu\leq\cf(\lambda)\leq\kappa$.\\
By Lemma \ref{lem3.16} we know that 
\[\cov(\lambda,\kappa^+,\kappa^+,\mu)\leq u(\kappa^+,\lambda)\leq
\lambda^+\]
and 
\[\lambda^+\leq u(\kappa^+,\lambda)\leq \big(\cov(\lambda,\kappa^+,
\kappa^+,\mu)\big)^{<\mu}.\]
Since $\lambda^{<\mu}=\lambda$, it follows that $\cov(\lambda,\kappa^+,
\kappa^+,\mu)=\lambda^+$. 
\medskip

\noindent{\sc Case:}\quad $\cf(\lambda)<\mu$.\\
By Lemmas \ref{lem3.13}, \ref{lem3.14} and \ref{lem3.16} we have 
\[\lambda\leq\cov(\lambda,\kappa^+,\kappa^+,\mu)\leq\sup\{u(\kappa^+,
\nu):\kappa<\nu<\lambda\}\leq\lambda,\]
and consequently $\cov(\lambda,\kappa^+,\kappa^+,\mu)=\lambda$.
\end{proof}

Thus, if the Singular Cardinals Hypothesis holds and $\lambda\geq 2^\kappa$,
then 
\[\cof^{<\mu}(\NSkl)<\cof(\NSkl)\quad\mbox{ if and only if }\quad \cf(
\lambda)<\mu\leq\kappa.\]  

\section{$\cov(\Mkk)$} 
Let us recall some basic facts and definitions related to the combinatorics
of the $\kappa$--meager ideal $\Mkk$ on $\bairek$. 

\begin{definition}
\label{d0.1}
\begin{enumerate}
\item {\em The Baire number\/} $n(X)$ of a topological space $X$ (also
called {\em the Novak number of $X$}) is the least number of nowhere dense 
subsets of $X$ needed to cover $X$.
\item For a topological space $X$ and a cardinal $\mu$, the $\mu$--complete 
ideal of subsets of $X$ generated by nowhere dense subsets of $X$ is denoted
by ${\bf M}_{<\mu}(X)$; ${\bf M}_{<\mu^+}(X)$ will be also denoted by ${\bf
M}_{\mu}(X)$. The ideal ${\bf M}_\mu(X)$ is {\em the ideal of $\mu$--meager
subsets of $X$}.  
\item The space $\bairek$ (respectively $\cantork$) is endowed with the
topology obtained by taking as basic open sets $\emptyset$ and $O_s$ for
$s\in\seqk$ (respectively $s\in {}^{\textstyle {<}\kappa}2$), where
$O_s=\{f\in \bairek: s\subseteq f\}$ (respectively $O_s=\{f\in \cantork:
s\subseteq f\}$). 
\item The ideals of $\kappa$--meager subsets of $\bairek$, $\cantork$ are
denoted by $\Mkkk$ and $\Mkk$, respectively.
\end{enumerate}
\end{definition}

\begin{remark}
\label{rem4.2}
\begin{enumerate}
\item Clearly, for a topological space $X$, $n(X)$ is the least number of
open dense subsets of $X$ with empty intersection. If $\mu\leq n(X)$, then
${\bf M}_{<\mu}(X)$ is a proper ideal (i.e., $X\notin{\bf M}_{<\mu}(X)$). 
\item Following the tradition of Set Theory of the Reals, we may consider
the covering number $\cov({\bf M}_{<\mu}(X))$ of the ideal ${\bf M}_{<\mu}
(X)$: 
\[\cov({\bf M}_{<\mu}(X))=\min\{|\cA|:\cA\subseteq {\bf M}_{<\mu}(X)\ \&\
\bigcup\cA=X\}.\]
By the definition, $n(X)=\cov({\bf M}_{<\aleph_0}(X))$. But also for every
$\mu<n(X)$ we have 
\[\cov({\bf M}_{<\mu}(X))=n(X)\quad\mbox{ and }\quad\cov({\bf M}_{<n(X)}(X))=
\cf(n(X)).\]   
\item Plainly, $n(\bairek)>\kappa$ and $n(\cantork)>\kappa$ (remember,
$\kappa$ is assumed to be regular).
\end{enumerate}
\end{remark}
 
\begin{lemma}
\label{lem4.3}
Suppose that $X$ is a topological space, $\mu<n(X)$, and $Y_\alpha$ are open
subsets of $X$ (for $\alpha<\mu$). Assume also that
$Y=\bigcap\limits_{\alpha<\mu} Y_\alpha$ is dense in $X$. Then, if $Y$ is
equipped with the subspace topology, $n(Y)=n(X)$. 
\end{lemma}

\begin{proof}
Let $U_\beta$ (for $\beta<n(X)$) be open dense subsets of $X$ such that
$\bigcap\limits_{\beta<n(X)}U_\beta=\emptyset$. Then $U_\beta\cap Y$ are
open dense subsets of $Y$ (remember $Y$ is dense) and $\bigcap\limits_{\beta
<n(X)}(U_\beta\cap Y)=\emptyset$. This shows that $n(Y)\leq n(X)$. 

Now, let $V_\beta\subseteq Y$ (for $\beta<n(Y)$) be open dense subsets of
$Y$ such that $\bigcap\limits_{\beta<n(Y)} V_\beta=\emptyset$. Take open
subsets $U_\beta$ of $X$ such that $V_\beta=U_\beta\cap Y$ -- clearly
$U_\beta$'s are dense in $X$ (as $Y$ is so). Then
$\emptyset=\bigcap\limits_{\beta<n(Y)}(U_\beta\cap Y)= 
\bigcap\limits_{\beta<n(Y)} U_\beta\cap\bigcap\limits_{\alpha<\mu}
Y_\alpha$, and hence $n(X)\leq n(Y) + \mu$ and therefore $n(X)\leq n(Y)$. 
\end{proof}

\begin{proposition}
\label{prop4.4}
$\cov(\Mkkk)=n(\bairek)=n(\cantork)=\cov({\bf M}_{\kappa,\kappa})$.
\end{proposition}

\begin{proof}
For $s\in {}^{\textstyle{<}\kappa}2$ and $\alpha<\kappa$ let $F(s,\alpha)\in
{}^{\textstyle{<}\kappa}2$ be such that $\lh(F(s,\alpha))=\lh(s)+\alpha+1$
and  
\[F(s,\alpha)\rest\lh(s)=s,\quad F(s,\alpha)\rest [\lh(s),\lh(s)+\alpha)
\equiv 1,\quad\mbox{and}\quad F(s,\alpha)(\lh(s)+\alpha)=0.\]
Now, let $\pi:\seqk\longrightarrow {}^{\textstyle{<}\kappa}2$ be such that 
\begin{itemize}
\item $\pi(\langle\rangle)=\langle\rangle$, $\pi(s\conc\langle\alpha\rangle)
=F(\pi(s),\alpha)$ for $s\in \seqk$, and
\item if $\langle s_\zeta:\zeta<\xi\rangle\subseteq\seqk$ is
$\vartriangleleft$--increasing, $\xi<\kappa$, $s=\bigcup\limits_{\zeta<\xi}
s_\zeta$,\\
then $\pi(s)=\bigcup\limits_{\zeta<\xi}\pi(s_\zeta)$. 
\end{itemize}
Then $\pi$ induces a mapping
\[\pi^*:\bairek\longrightarrow\cantork:\eta\mapsto\bigcup_{\zeta<\kappa}
\pi(\eta\rest\zeta).\] 
The range of $\pi^*$ is  
\[\Rng(\pi^*)=\{\rho\in\cantork: (\forall\alpha<\kappa)(\exists\beta<
\kappa)(\alpha<\beta\ \&\ \rho(\beta)=0)\}.\]
Plainly, $\Rng(\pi^*)$ is the intersection of $\kappa$ many open dense
subsets of $\cantork$. Moreover, $\pi^*$ is a homeomorphism from
$\bairek$ onto $\Rng(\pi^*)$. Therefore, using Lemma \ref{lem4.3}, we get
$n(\bairek)=n(\Rng(\pi^*))=n(\cantork)$. The rest should be clear (remember
Remark \ref{rem4.2}(2,3)). 
\end{proof}

\begin{proposition}
\label{trivial}
$\cov(\Mkk)\leq\gd_\kappa$.
\end{proposition}

\begin{definition}
\label{d0.2}
$\Cmk$ is the forcing notion for adding $\mu$ Cohen functions in $\bairek$
with ${<}\kappa$--support. Thus a condition in $\Cmk$ is a function $q$ such
that   
\[\Dom(q)\subseteq\mu\times\kappa,\qquad \Rng(q)\subseteq\kappa\quad
\mbox{and}\quad |q|<\kappa.\]
The order of $\Cmk$ is the inclusion.
\end{definition}

\begin{proposition}
\label{p0.3}
Assume $2^{<\kappa}=\kappa<\mu$. Then $\forces_{\Cmk}$`` $\cov(\Mkk)\geq\mu$
''.
\end{proposition}

\section{Manageable forcing notions}
In this section we introduce a property of forcing notions which is crucial
for the consistency results presented later: $(\theta,\mu,
\kappa)$--manageability. This property has three ingredients: an iterable
variant of $\kappa^+$--cc (see Definition \ref{288cc}),
$\kappa$--completeness and a special property implying preservation of
non-meagerness of subsets of $\bairek$ (see Proposition \ref{p1.5}). Since
later we will work with ${<}\kappa$--support iterations, we also prove a
suitable preservation theorem (see Theorem \ref{t1.6}).  

From now on we will always assume that our fixed (uncountable) regular
cardinal $\kappa$ satisfies $2^{<\kappa}=\kappa$ (so also $\kappa^{<\kappa}
=\kappa$). 

\begin{definition}
[See Shelah {\cite[Definition 1.1]{Sh:288}} and {\cite[Definition
7]{Sh:546}}] 
\label{288cc}
Let $\bP$ be a forcing notion, and $\varepsilon<\kappa$ be a limit ordinal. 
\begin{enumerate}
\item We define a game $\Gcc(\bP)$ of two players, Player I and Player II. A
play lasts $\varepsilon$ steps, and at each stage $\alpha<\varepsilon$ of
the play $\bar{q}^\alpha,\bar{p}^\alpha,\varphi^\alpha$ are chosen so that:
\begin{itemize}
\item $\bar{q}^0=\langle \emptyset_\bP:i<\kappa^+\rangle$, $\varphi^0:
\kappa^+\longrightarrow\kappa^+:i\mapsto 0$;
\item If $\alpha>0$, then Player I picks $\bar{q}^\alpha,\varphi^\alpha$
such that   
\begin{enumerate}
\item[(i)] $\bar{q}^\alpha=\langle q^\alpha_i:i<\kappa^+\rangle\subseteq
\bP$ satisfies
\[(\forall \beta<\alpha)(\forall i<\kappa^+)(p^\beta_i\leq q^\alpha_i),\]  
\item[(ii)] $\varphi^\alpha:\kappa^+\longrightarrow\kappa^+$ is regressive,
i.e., $(\forall i<\kappa^+)(\varphi^\alpha(i)<1+i)$; 
\end{enumerate}
\item Player II answers choosing a sequence $\bar{p}^\alpha=\langle
p^\alpha_i: i<\kappa^+\rangle\subseteq\bP$ such that $(\forall i<\kappa^+)
(q^\alpha_i\leq p^\alpha_i)$.
\end{itemize}
If at some stage of the game Player I does not have any legal move, then he
looses. If the game lasted $\varepsilon$ steps, Player I wins a play
$\langle\bar{q}^\alpha,\bar{p}^\alpha,\varphi^\alpha:\alpha<\varepsilon
\rangle$ if there is a club $C$ of $\kappa^+$ such that for each distinct
members $i,j$ of $C$ satisfying $\cf(i)=\cf(j)=\kappa$ and
$(\forall\alpha<\varepsilon)(\varphi^\alpha(i)=\varphi^\alpha(j))$, the set
\[\{p^\alpha_i:\alpha<\varepsilon\}\cup \{p^\alpha_j: \alpha<\varepsilon\}\]
has an upper bound in $\bP$. 

\item The forcing notion $\bP$ satisfies condition $(*)^\varepsilon_\kappa$
if Player I has a winning strategy in the game $\Gcc(\bP)$.
\end{enumerate}
\end{definition}

\begin{remark}
Condition $(*)^\varepsilon_\kappa$ is a strong version of $\kappa^+$--cc
(easily, if $\varepsilon<\kappa$ is limit, $\kappa^\varepsilon=\kappa$, and
$\bP$ satisfies $(*)^\varepsilon_\kappa$, then $\bP$ satisfies
$\kappa^+$--cc). This condition was used in a number of papers, e.g., to
obtain a series of consistency results on partition relations; see Shelah
and Stanley \cite{ShSt:154}, \cite{ShSt:154a}, Shelah \cite{Sh:80},
\cite{Sh:288}, \cite{Sh:546}. Its primary use comes from the fact that it is  
preserved in ${<}\kappa$--support iterations.
\end{remark}

\begin{proposition}
[See Shelah {\cite[Iteration Lemma 1.3]{Sh:288}} and {\cite[Theorem
35]{Sh:546}}] 
\label{288it}
Let $\varepsilon<\kappa$ be a limit ordinal,
$\kappa=\kappa^{<\kappa}$. Suppose that $\bar{\bQ}=\langle\bP_\xi,\dbQ_\xi:
\xi<\gamma\rangle$ is a ${<}\kappa$--support iteration such that for each
$\xi<\gamma$  
\[\forces_{\bP_\xi}\mbox{`` }\dbQ_\xi\mbox{ satisfies
}(*)^\varepsilon_\kappa\mbox{ ''.}\]
Then $\bP_\gamma$ satisfies $(*)^\varepsilon_\kappa$.
\end{proposition}

\begin{definition}
\label{d1.1}
A forcing notion $\bP$ is {\em $\theta$--complete\/} if every
$\leq_{\bP}$--increasing chain of length less than $\theta$ has an upper
bound in $\bP$. It is {\em $\theta$--lub--complete\/} if every
$\leq_{\bP}$--increasing chain of length less than $\theta$ has a least
upper bound in $\bP$.
\end{definition}

\begin{definition}
\label{d1.2}
Let $\theta$ and $\mu$ be cardinals such that $\theta<\kappa=2^{<\kappa}$
and $\mu^{<\kappa}=\mu$. Let $\bP$ be a $\theta^+$--lub--complete forcing
notion.   
\begin{enumerate}
\item A model $N\prec (\cH(\chi),\in,<^*_\chi)$ is {\em $(\bP,\kappa,
\mu)$--relevant} if $\bP,\mu\in N$, $\mu\subseteq N$, $|N|=\mu$ and
${}^{\textstyle {<}\kappa}N\subseteq N$. 
\item For a $(\bP,\kappa,\mu)$--relevant model $N$ we define a game
$\Game^m(N,\theta,\bP)$ of two players, He and She, as follows. A play
lasts $\theta$ moves, and in the $i^{\rm th}$ move conditions
$p_i,q_i\in\bP$ are chosen so that:
\begin{itemize}
\item $q_i\in N\cap\bP$, $q_i\leq p_i$,
\item $(\forall j<i)(q_j\leq q_i\ \&\ p_j\leq p_i)$,
\item She chooses $p_i,q_i$ if $i$ is odd, He picks $p_i,q_i$ if $i$ is
even. 
\end{itemize}
She wins the play $\langle q_i,p_i:i<\theta\rangle$ whenever

{\em if}\quad $p_\theta$ is a least upper bound of $\langle
p_i:i<\theta\rangle$, and $q_\theta$ is a least upper bound of $\langle
q_i:i<\theta\rangle$, 

{\em then}\ \ $(\forall q\in N\cap\bP)(q_\theta\leq q\ \Rightarrow\
q,p_\theta\mbox{ are compatible })$.
\item The forcing notion $\bP$ is {\em weakly $(\theta,\mu,
\kappa)$--manageable\/} if (it is $\theta^+$--lub--complete and) there is an
$x\in\cH(\chi)$ (called {\em a witness\/}) such that for every
$(\bP,\kappa,\mu)$--relevant model $N\prec\cH(\chi)$ with $x\in N$, She has
a winning strategy in the game $\Game^m(N,\theta,\bP)$. 
\item  The forcing notion $\bP$ is {\em $(\theta,\mu,\kappa)$--manageable
\/} if it is $\kappa$--complete, weakly $(\theta,\mu,\kappa)$--manageable,
and satisfies the condition $(*)^\theta_\kappa$. 
\end{enumerate}
\end{definition}

\begin{remark}
Suppose that $\bP$ is $\theta^+$--lub--complete and $N$ is
$(P,\kappa,\mu)$--relevant. Then both players have always legal moves in
the game $\Game^m(N,\theta,\bP)$. Moreover, if $\langle q_i,p_i:i<\theta
\rangle$ is a (legal) play of $\Game^m(N,\theta,\bP)$, then there are least
upper bounds $q_\theta\in N\cap\bP$ of $\langle q_i:i<\theta\rangle$, and
$p_\theta\in\bP$ of $\langle p_i:i<\theta\rangle$ (and $q_\theta\leq
p_\theta$). 
\end{remark}

\begin{definition}
\label{d1.3A}
Let $N$ be a $(P,\kappa,\mu)$--relevant model, and let $q\in N\cap\bP$,
$p\in \bP$ be such that $q\leq p$. We say that a pair $(q^*,p^*)$ is {\em an   
$N$--cover for $(q,p)$}, if 
\begin{itemize}
\item $q\leq q^*\in N\cap\bP$, $p\leq p^*\in\bP$, $q^*\leq p^*$, and
\item every condition $q'\in N\cap\bP$ stronger than $q^*$ is compatible with
$p^*$. 
\end{itemize}
\end{definition}

\begin{lemma}
\label{p1.4}
Suppose that $\bP$ is a $\theta^+$--lub--complete forcing notion, $N$ is a
$(P,\kappa,\mu)$--relevant model, and She has a winning strategy in the
game $\Game^m(N,\theta,\bP)$. Then:
\begin{enumerate}
\item For all conditions $q\in N\cap\bP$ and $p\in\bP$ such that $q\leq p$,
there is an $N$--cover $(q^*,p^*)$ for $(q,p)$. 
\item $N\cap\bP\lesdot \bP$.
\end{enumerate}
\end{lemma}

\begin{proof}
1)\quad Consider a play $\langle q_i,p_i:i<\theta\rangle$ of $\Game^m(N,
\theta,\bP)$ in which He starts with $q_0=q$, $p_0=p$, and then he always
plays the $<^*_\chi$--first legal moves, and She uses her winning strategy. 
Let $q^*\in N\cap\bP$, $p^*\in\bP$ be least upper bounds of $\langle q_i:i<
\theta\rangle$, $\langle p_i:i<\theta\rangle$, respectively. Plainly, as
She won the play, the pair $(q^*,p^*)$ is an $N$--cover for $(q,p)$. 

\noindent 2)\quad Suppose that $\cA\subseteq N\cap\bP$ is a maximal
antichain in $N\cap\bP$, but $p\in \bP$ is incompatible with all members of
$\cA$. Let $(q^*,p^*)$ be an $N$--cover for $(\emptyset_\bP,p)$. The
condition $q^*$ is compatible with some $q\in\cA$, so let $q^+\in N\cap \bP$
be such that $q^+\geq q^*$, $q^+\geq q\in\cA$. By the choice of $(q^*,p^*)$
we know that the conditions $q^+$ and $p^*$ are compatible, and hence $q$
and $p$ are compatible. A contradiction. 

The rest follows from the elementarity of $N$.
\end{proof}

\begin{proposition}
\label{p1.5}
Assume $\theta<\kappa=2^{<\kappa}\leq\mu=\mu^{<\kappa}<\tau$. Suppose that
a set $Y\subseteq\bairek$ cannot be covered by the union of less than $\tau$
nowhere dense subsets of $\bairek$, and $\bP$ is a weakly $(\theta,\mu,
\kappa)$--manageable forcing notion not collapsing cardinals. Then 
\[\forces_{\bP}\mbox{`` $Y$ is not the union of $<\tau$ nowhere dense
subsets of $\bairek$ ''.}\]
\end{proposition}

\begin{proof}
Let $\bP$ be weakly $(\theta,\mu,\kappa)$--manageable with a witness
$x\in\cH(\chi)$. Suppose toward contradiction that a condition $q\in\bP$ is
such that 
\[q\forces\mbox{`` $Y$ is the union of $<\tau$ nowhere dense subsets of
$\bairek$ ''}.\]  
Passing to a stronger condition if needed, we may assume that for some
$\iota<\tau$ and $\bP$--names $\dot{A}_\xi$ (for $\xi<\iota$) we have:
\[q\forces\mbox{`` }\dot{A}_\xi\subseteq {}^{\textstyle{<}\kappa}\kappa\ \&\ 
(\forall s\in {}^{\textstyle{<}\kappa}\kappa)(\exists t\in\dot{A}_\xi)(s
\subseteq t)\mbox{ ''}\]
and 
\[q\forces\mbox{`` }(\forall y\in Y)(\exists\xi<\iota)(\forall t\in
\dot{A}_\xi)(t\nsubseteq y)\mbox{ ''.}\]
For each $\zeta<\iota$ pick a $(\bP,\kappa,\mu)$--relevant model
$N_\zeta\prec (\cH(\chi),\in,<^*_\chi)$ such that $q,\langle\dot{A}_\xi:
\xi<\kappa\rangle,x,\zeta\in N_\zeta$. Then $|\bigcup\limits_{\zeta<
\iota} N_\zeta|=\iota\cdot\mu<\tau$, so we may pick a $y\in Y$ 
such that $y\in\cO$ for all open dense subsets $\cO$ of $\bairek$ from 
$\bigcup\limits_{\zeta<\iota}N_\zeta$. By our assumptions, there are
$\xi<\iota$ and $p\geq q$ such that 
\[p\forces\mbox{`` }(\forall t\in \dot{A}_\xi)(t\nsubseteq y)\mbox{ ''.}\]
Let $(q^*,p^*)$ be an $N_\xi$--cover for $(q,p)$ (there is one by
Lemma \ref{p1.4}(1)). Put
\[A=\{s\in{}^{\textstyle {<}\kappa}\kappa:(\exists q'\geq q^*)(q'\forces
s\in\dot{A}_\xi)\}.\] 
Clearly $A\in N_\xi$, $A\subseteq N_\xi$, and $\cO=\bigcup\limits_{s\in A}
O_s\in N_\xi$ is an open dense subset of $\bairek$. Hence $s\subseteq y$ for
some $s\in A$. Let $q'\in N_\xi\cap\bP$ be a condition stronger than $q^*$
and such that $q'\forces s\in\dot{A}_\xi$. The condition $q'$ is compatible
with  $p^*$, and so with $p$. Take a condition $q^+$ stronger than both $q'$
and $p$. Then
\[q^+\forces\mbox{`` }s\in\dot{A}_\xi\ \&\ s\subseteq y\mbox{ ''}
\quad\mbox{ and }\quad q^+\forces\mbox{`` }(\forall t\in\dot{A}_\xi)(t
\nsubseteq y)\mbox{ '',}\]
a contradiction.
\end{proof}

\begin{corollary}
\label{addedc}
Suppose that $\theta<\kappa=2^{<\kappa}\leq\mu=\mu^{<\kappa}$ and
$\cov(\Mkk)>\mu$. Let $\bP$ be a $(\theta,\mu,\kappa)$--manageable forcing
notion. Then 
\[\forces_{\bP}\mbox{`` }\big(\cov(\Mkk)\big)^\bV\leq \cov(\Mkk)\mbox{
''.}\] 
\end{corollary}

\begin{proof}
Rememebering Proposition \ref{prop4.4}, apply Proposition \ref{p1.5} to
$\tau=\cov(\Mkk)=\cov(\Mkkk)$ and $Y=\bairek$ to get 
\[\forces_{\bP}\mbox{`` $(\bairek)^\bV$ is not the union of $<\tau$ nowhere
dense sets ''.}\]
But this clearly implies $\forces_{\bP}\mbox{`` }\tau\leq
\cov(\Mkk)=\cov(\Mkkk)\mbox{ ''}$. 
\end{proof}

\begin{theorem}
\label{t1.6}
Assume that $\theta<\kappa=2^{<\kappa}\leq\mu=\mu^{<\kappa}$. Let
$\bar{\bQ}=\langle\bP_\xi,\dbQ_\xi:\xi<\gamma\rangle$ be
${<}\kappa$--support iteration such that for each $\xi<\gamma$
\[\forces_{\bP_\xi}\mbox{`` $\dbQ_\xi$ is $(\theta,\mu,\kappa)$--manageable  
''}.\]
Then $\bP_\gamma$ is $(\theta,\mu,\kappa)$--manageable.
\end{theorem}

\begin{proof}
Let $\theta,\kappa,\mu$ and $\bar{\bQ}$ be as in the assumptions of the
theorem. 

First note that ${<}\kappa$--support iterations of $\kappa$--complete
forcing notions satisfying the condition $(*)^\theta_\kappa$ are
$\kappa$--complete $\kappa^+$--cc (as $\kappa^{<\kappa}=\kappa$; remember
Proposition \ref{288it}).  Therefore no such iteration collapses cardinals
nor changes cofinalities nor adds sequences of ordinals of length
$<\kappa$. Hence the assumed properties of $\theta,\kappa$ and $\mu$ hold in
all intermediate extensions $\bV^{\bP_\xi}$ and our assumption on
$\dbQ_\xi$'s is meaningful.  

Plainly, $\bP_\gamma$ is $\kappa$--complete, $\theta^+$--lub--complete and
satisfies condition $(*)^\theta_\kappa$. We have to show that $\bP_\gamma$
is weakly $(\theta,\mu,\kappa)$--manageable.

For $\xi<\gamma$ let $\dot{x}_\xi$ be a $\bP_\xi$--name for a witness for
$\dbQ_\xi$ being weakly manageable and let $\bar{x}=\langle\dot{x}_\xi:\xi
<\gamma\rangle$. Suppose that $N\prec (\cH(\chi),\in,<^*_\chi)$ is a
$(\bP_\gamma,\kappa,\mu)$--relevant model such that $(\bar{x},\bar{\bQ})\in
N$. 

Since $\bP_\xi$ satisfies $\kappa^+$--cc (and $\kappa+1\subseteq N$) we know 
that if $\xi\in N\cap \gamma$ and $G_\xi\subseteq \bP_\xi$ is generic over
$\bV$, then in $\bV[G_\xi]$ we have:
\[N[G_\xi]\cap\bV=N\quad\mbox{ and }\quad N[G_\xi]\prec (\cH(\chi),\in,
<^*_\chi)^{\bV[G_\xi]}\quad\mbox{ and }\quad {}^{\textstyle {<}\kappa}N[
G_\xi]\subseteq N[G_\xi]\]
(remember that $\bP_\xi$ is $\kappa$--complete). As clearly $\dbQ_\xi^{
G_\xi}\in N[G_\xi]$, we conclude that $N[G_\xi]$ is $(\dbQ_\xi^{G_\xi},
\kappa,\mu)$--relevant, and $\dot{x}_\xi^{G_\xi}\in N[G_\xi]$. Therefore,
She has a winning strategy in the game $\Game^m(N[G_\xi],\theta,
\dbQ_\xi^{G_\xi})$. Let $\nst_\xi$ be a $\bP_\xi$--name for such a
strategy. We may assume that the strategy $\nst_\xi$ is such that 
\begin{enumerate}
\item[$(*)$] if $i<\theta$ is even and $q_i=p_i=\dot{\emptyset}_{
\dbQ_\xi}$,\\
then $\nst_\xi$ instructs Her to play $q_{i+1}=p_{i+1}=
\dot{\emptyset}_{\dbQ_\xi}$. 
\end{enumerate}
We define a strategy $\st$ for Her in the game $\Game^m(N,\theta,
\bP_\gamma)$ as follows. At an odd stage $i<\theta$ of the game, the
strategy $\st$ first instructs Her to choose (side) conditions $q_i^-,p_i^-
\in\bP_\gamma$ and only then pick conditions $q_i\in N\cap\bP_\gamma$ and
$p_i\in\bP_\gamma$ which are to be played. These conditions will be chosen
so that if $\langle q_j,p_j:j<i\rangle$ is a legal play of $\Game^m(N,
\theta,\bP_\gamma)$ in which She uses $\st$, and $q_j^-,p_j^-$ are the side 
conditions picked by her (for odd $j\leq i$), then
\begin{enumerate}
\item[$(\alpha)_i$] $\Dom(q_i)=\Dom(q_i^-)=\Dom(p_{i-1})\cap N$,
$\Dom(p_i^-)=\Dom(p_{i-1})$,  
\item[$(\beta)_i$]  $p_{i-1}\leq p_i^-\leq p_i$, $q_{i-1}\leq q_i^-\leq
p_i^-$, $q_{i-1}\leq q_i\leq p_i$,
\item[$(\gamma)_i$] letting $(q_j^*,p_j^*)$ be $(q_j^-,p_j^-)$ if $j\leq i$
is odd and $(q_j,p_j)$ if $j<i$ is even, for every $\xi\in\Dom(q_i)$ we have  
\[\begin{array}{ll}
p_i\rest\xi\forces_{\bP_\xi}&\mbox{`` }q_i(\xi)=q_i^-(\xi)\mbox{ and}\\
&\ \ \ \mbox{the sequence }\langle q_j^*(\xi),p_j^*(\xi):j\leq i\rangle
\mbox{ is a legal play of}\\
&\ \ \ \Game^m(N[\Gamma_{\bP_\xi}],\theta,\dbQ_\xi)\mbox{ in which She uses
the strategy $\nst_\xi$ ''.}
  \end{array}\]
\end{enumerate}
So suppose that $i<\theta$ is odd, $\langle q_j,p_j:j<i\rangle$ is a
partial play of $\Game^m(N,\theta,\bP_\gamma)$ in which She uses $\st$ (and
the side conditions for odd $j<i$ are $q^-_j,p^-_j$), and the clauses 
$(\alpha)_j$--$(\gamma)_j$ hold for all odd $j<i$. Let $(q_j^*,p_j^*)$ be
$(q_j^-,p_j^-)$ if $j<i$ is odd and $(q_j,p_j)$ if $j<i$ is even. 

We first declare that $\Dom(q^-_i)=\Dom(p_{i-1})\cap N$, $\Dom(p_i^-)=\Dom(
p_{i-1})$ and $p_i^-(\zeta)=p_{i-1}(\zeta)$ for all $\zeta\in\Dom(p_i^-)
\setminus N$. Next, by induction on $\xi\in \Dom(q^-_i)$ we define
$q^-_i(\xi),p^-_i(\xi)$. So suppose that $\xi\in\Dom(q^-_i)$ and $q_i^-\rest
\xi, p_i^-\rest\xi$ have been defined so that $q_{i-1}\rest\xi\leq q_i^-
\rest\xi$, $p_{i-1}\rest\xi\leq p_i^-\rest\xi$. Then, by clauses
$(\gamma)_j$,   
\[\begin{array}{ll}
p_i^-\rest\xi\forces_{\bP_\xi}&\mbox{`` the sequence }\langle q_j^*(\xi),
p_j^*(\xi):j<i\rangle\mbox{ is a legal play of}\\
&\ \ \ \Game^m(N[\Gamma_{\bP_\xi}],\theta,\dbQ_\xi)\mbox{ in which She uses
the strategy $\nst_\xi$ ''.}
  \end{array}\]
(Remember our assumption $(*)$ on $\nst_\xi$ and our convention regarding
$\emptyset_\bP$ stated in Notation \ref{notation}(1).) 
Let $q_i^-(\xi)$ and $p_i^-(\xi)$ be $\bP_\xi$--names for members of
$\dbQ_\xi$ such that
\[q_i^-\rest\xi\forces_{\bP_\xi}\mbox{`` }q_i^-(\xi)\in N[\Gamma_{\bP_\xi}]\
\&\ q_{i-1}(\xi)\leq q^-_i(\xi)\mbox{ '',}\]
and 
\[\begin{array}{ll}
p_i^-\rest\xi\forces_{\bP_\xi}&\mbox{`` }(q_i^-(\xi),p_i^-(\xi))\mbox{ is
what $\nst_\xi$ says Her to play}\\
&\ \mbox{ as the answer to }\langle q_j^*(\xi),p_j^*(\xi):j<i\rangle
\mbox{ ''.} 
  \end{array}\]
(So $q_i^-(\xi)$ is a name for a member of $N[\Gamma_{\bP_\xi}]$, but it
does not have to be from $N$.) This completes the definition of $q_i^-,p_i^-
\in\bP_\gamma$. Now we use the fact that $\bP_\gamma$ is $\kappa$--complete
and $|\Dom(q^-_i)|<\kappa$ to pick a condition $p_i\in\bP_\gamma$ stronger
than $p_i^-$ and names $\dot{\tau}_\xi\in N$ (for $\xi\in\Dom(q^-_i)$) such
that $p_i\rest\xi\forces_{\bP_\xi}$`` $q_i^-(\xi)=\dot{\tau}_\xi$ ''. Since
${}^{\textstyle {<}\kappa}N\subseteq N$, the sequence $\langle\dot{\tau}_\xi:
\xi\in\Dom(q^-_i)\rangle$ is in $N$. Hence we may find a condition $q_i\in
N\cap\bP_\gamma$ such that
\begin{itemize}
\item $\Dom(q_i)=\Dom(q_i^-)$, and 
\item for each $\xi\in\Dom(q_i)$, 
\[\forces_{\bP_\xi}\mbox{`` if }\dot{\tau}_\xi\geq q_{i-1}(\xi),\mbox{ then
}q_i(\xi)=\dot{\tau}_\xi,\mbox{ otherwise }q_i(\xi)=q_{i-1}(\xi)\mbox{ ''.}\]
\end{itemize}
(For definitiveness we pick the $<^*_\chi$--first $p_i,q_i$ as above.) It
should be clear that $q_i^-,q_i,p_i^-,p_i$ satisfy conditions
$(\alpha)_i$--$(\gamma)_i$. This finishes the description of the strategy
$\st$. Let us argue that it is a winning strategy for Her.

To this end suppose that $\langle q_i,p_i:i<\theta\rangle$ is the result of
a play of $\Game^m(N,\theta,\bP_\gamma)$ in which She uses $\st$. Let
$q_\theta,p_\theta\in \bP_\gamma$ be least upper bounds of $\langle q_i:
i<\theta\rangle$, $\langle p_i:i<\theta\rangle$, respectively. Then for
every $\xi\in\Dom(p_\theta)$ we have 
\[p_\theta\rest\xi\forces_{\bP_\xi}\mbox{`` }p_\theta(\xi)\mbox{ is a
least upper bound of }\langle p_i(\xi):i<\theta\rangle\mbox{ ''.}\]
We may also assume that $\Dom(q_\theta)=\bigcup\limits_{i<\theta}\Dom(q_i)
=\Dom(p_\theta)\cap N$.

Let $q\in N\cap\bP_\gamma$ be a condition stronger than $q_\theta$ (and
thus stronger than all $q_i$ for $i<\theta$). We define a condition
$p\in\bP_\gamma$ as follows. First, we declare that $\Dom(p)=\Dom(q)\cup
\Dom(p_\theta)$, and $p(\xi)=q(\xi)$ for $\xi\in\Dom(q)\setminus \Dom(
p_\theta)$, and $p(\xi)=p_\theta(\xi)$ for $\xi\in\Dom(p_\theta)\setminus
\Dom(q)=\Dom(p_\theta)\setminus N$. Now suppose that $\xi\in\Dom(q)\cap
\Dom(p_\theta)$ and we have already defined $p\rest\xi$ so that $q\rest\xi
\leq p\rest\xi$ and $p_\theta\rest\xi\leq p\rest\xi$. Then, by our choices,
\[\begin{array}{ll}
p\rest\xi\forces_{\bP_\xi}&\mbox{`` the sequence }\langle q_j(\xi),
p_j^*(\xi):j<\theta\rangle\mbox{ is a legal play of}\\
&\ \ \ \Game^m(N[\Gamma_{\bP_\xi}],\theta,\dbQ_\xi)\mbox{ in which She uses 
the strategy $\nst_\xi$, and}\\
&\ \ \ q(\xi)\in N[\Gamma_{\bP_\xi}]\mbox{ is stronger than all }q_j(\xi)
\mbox{ for }j<\theta\mbox{ ''.} 
  \end{array}\]
(Above, $p^*_j$ are as in the definition of $\st$: either $p_j$ or $p_j^-$,
depending on the parity of $j$.) Consequently, 
\[p\rest\xi\forces\mbox{`` }q(\xi)\mbox{ and }p_\theta(\xi)\mbox{ are
compatible '',}\]
so we may pick a $\bP_\xi$--name $p(\xi)$ for a condition in $\dbQ_\xi$ such
that 
\[p\rest\xi\forces\mbox{`` }q(\xi)\leq p(\xi)\mbox{ and }p_\theta(\xi)\leq
p(\xi)\mbox{ ''.}\] 
This completes the choice of $p\in\bP_\gamma$. Plainly, $p$ is an upper
bound of $q$ and $p_\theta$ showing that they are compatible.
\end{proof}

\section{The one-step forcing}
In this section we introduce a forcing notion $\bQ$ for adding a small
family of functions in $\bairek$ which $\tau$--dominates
$\bairek\cap\bV$. Iterating this type of forcing notions we will get models
with $\gd_\kappa^\tau$ small. Our forcing is (of course) manageable for
suitable parameters, and thus it preserves non-meagerness of subsets of
$\kappa$. Throughout this section we assume the following.

\begin{context}
\label{c3.1}
\begin{enumerate}
\item[(i)]   $\kappa=2^{<\kappa}$, $\tau=\cf(\tau)<\kappa$, 
\item[(ii)]  $\bar{\mu}=\langle\mu_\alpha:\alpha<\tau\rangle$ is an
increasing sequence of regular cardinals, $\kappa\leq \mu_0$, 
\item[(iii)] $|\prod\limits_{\alpha<\tau}\mu_\alpha|=2^\kappa$ and
$\pi:\prod\limits_{\alpha<\tau}\mu_\alpha\longrightarrow{}^{\textstyle
\kappa}\kappa$ is a bijection.
\end{enumerate}
We will write $\pi_\eta$ for $\pi(\eta)$. Also for a set $u\subseteq
\prod\limits_{\alpha<\tau}\mu_\alpha$ we let 
\[T(u)\stackrel{\rm def}{=}\{\eta\rest\alpha:\alpha<\tau\ \&\ \eta\in
u\}\]
\end{context}

\begin{definition}
\label{d3.2}
\begin{enumerate}
\item We define a forcing notion $\bQ=\bQ(\pi,\bar{\mu},\kappa)$ as
follows.\\ 
{\bf A condition} in $\bQ$ is a tuple $p=(i,u,\bar{f},g)=(i^p,u^p,
\bar{f}^p,g^p)$ such that 
\begin{enumerate}
\item[(a)] $i<\kappa$, $u\in \cP_\kappa\big(\prod\limits_{\alpha<\tau}
\mu_\alpha\big)$, 
\item[(b)] $\bar{f}=\langle f_\sigma:\sigma\in T(u)\rangle$ and $f_\sigma:
i\longrightarrow\kappa$ for $\sigma\in T(u)$,
\item[(c)] $g:u\longrightarrow i+1$, and if $\eta\in u$, $g(\eta)\leq j<i$,
then 
\[\pi_\eta(j)<\sup\{f_{\eta\rest\alpha}(j):\alpha<\tau\}.\]
\end{enumerate}
{\bf The order of $\bQ$} is such that for $p,q\in\bQ$ we have 

$p\leq q$ if and only if 

$i^p\leq i^q$, $u^p\subseteq u^q$, $g^p\subseteq g^q$ and $f^p_\sigma
\subseteq f^q_\sigma$ for $\sigma\in T(u^p)$.
\item For a set $U\subseteq \prod\limits_{\alpha<\tau}\mu_\alpha$ we
let $\bQ\rest U=\{p\in \bQ: u^p\subseteq U\}$, and for a condition $q\in\bQ$
we put 
\[q\rest U=\big(i^q,u^q\cap U,\bar{f}^q\rest T(u^q\cap U),g^q\rest (u^q\cap
U)\big).\]  
\end{enumerate}
\end{definition}

\begin{proposition}
\label{p3.3}
\begin{enumerate}
\item $\bQ$ is a $\kappa$--lub--complete forcing notion of size $2^\kappa$.
\item Let $U\subseteq \prod\limits_{\alpha<\tau}\mu_\alpha$ be of size
$\leq\kappa$. Then $|\bQ\rest U|\leq \kappa$.
\end{enumerate}
\end{proposition}

\begin{proof}
1)\quad Plainly, $(\bQ,\leq)$ is a partial order of size $2^\kappa$. To
prove the completeness suppose that $\langle p_\xi:\xi<\xi^*\rangle$ is an 
$\leq$--increasing sequence of members of $\bQ$ and $\xi^*<\kappa$. Put 
\[i^q=\sup_{\xi<\xi^*}i^{p_\xi},\quad u^q=\bigcup_{\xi<\xi^*}u^{p_\xi},
\quad g^q=\bigcup_{\xi<\xi^*} g^{p_\xi}\]
and $f^q_\sigma=\bigcup\{f^{p_\xi}_\sigma:\xi<\xi^*\ \&\ \sigma\in
T(u^{p_\xi})\}$ for $\sigma\in T(u^q)$. Clearly $q=(i^q,u^q,\bar{f}^q,g^q) 
\in\bQ$ is the least upper bound of $\langle p_\xi:\xi<\xi^*\rangle$. 

\noindent 2)\quad Should be clear.
\end{proof}

\begin{proposition}
\label{p633}
The forcing notion $\bQ$ satisfies the condition $(*)^\theta_\kappa$ (see
\ref{288cc}(2)) for any limit ordinal $\varepsilon<\kappa$. 
\end{proposition}

\begin{proof}
Let $\varepsilon<\kappa$ be a limit ordinal. To give the winning strategy
for Player I in the game $\Gcc(\bQ)$ we need two technical observations. 

\begin{claim}
\label{cl1}
If $p,q\in\bQ$ are such that $i^p=i^q$ and $g^p\rest (u^p\cap u^q)=g^q\rest
(u^p\cap u^q)$ and $f^p_\sigma=f^q_\sigma$ for $\sigma\in T(u^p)\cap
T(u^q)$, then the conditions $p,q$ have a least upper bound. 
\end{claim}

\begin{proof}[Proof of the Claim]
Let $i^r=i^p=i^q$, $u^r=u^p\cup u^q$, $g^r=g^p\cup g^q$ and 
\[f^r_\sigma=\left\{\begin{array}{ll}
f^p_\sigma &\mbox{ if } \sigma\in T(u^p)\\
f^q_\sigma &\mbox{ if } \sigma\in T(u^q).
		    \end{array}\right.\]
Then $r=(i^r,u^r,\bar{f}^r,g^r)\in\bQ$ is the least upper bound of $p,q$.
\end{proof}

\begin{claim}
\label{cl2}
Suppose $\bar{q}=\langle q_j:j<\kappa^+\rangle\subseteq\bQ$. Then there is a
regressive function $\varphi_{\bar{q}}:\kappa^+\longrightarrow\kappa^+$ such
that\\ 
if $j<j'<\kappa^+$, $\cf(j)=\cf(j')=\kappa$ and $\varphi_{\bar{q}}(j)=
\varphi_{\bar{q}}(j')$,\\ 
then $i^{q_j}=i^{q_{j'}}$, and $g^{q_j}\rest (u^{q_j}\cap u^{q_{j'}})=g^{
q_{j'}}\rest (u^{q_j}\cap u^{q_{j'}})$, and $f^{q_j}_\sigma=
f^{q_{j'}}_\sigma$ for $\sigma\in T(u^{q_j})\cap T(u^{q_{j'}})$.
\end{claim}

\begin{proof}[Proof of the Claim]
Take a sequence $\langle\eta_\xi:\xi<\kappa^+\rangle\subseteq
\prod\limits_{\alpha<\tau}\mu_\alpha$ such that for each $j<\kappa^+$
of cofinality $\kappa$ and an $\alpha<j$ we have $u^{q_\alpha}\subseteq
\{\eta_\xi:\xi<j\}$. Let $U=\{\eta_\xi:\xi<\kappa^+\}$ and $U_j=\{\eta_\xi:
\xi<j\}$ for $j<\kappa^+$. By \ref{p3.3}(2) we know that $|\bQ\rest U_j|\leq
\kappa$ (for $j<\kappa^+$)  and $|\bQ\rest U|\leq \kappa^+$, and hence we
may pick a mapping $\psi_0:\kappa^+\longrightarrow\bQ\rest U$ such that  
\[\big(\forall j<\kappa^+\big)\big(\cf(j)=\kappa\ \Rightarrow\ \Rng(\psi_0
\rest j)=\bQ\rest U_j \big).\] 
Also, for $j<\kappa^+$, let $F(U_j)$ be the set 
\[\big\{\bar{f}=\langle f_\sigma:\sigma\in\Dom(\bar{f})\rangle:\Dom(\bar{f})
\in \cP_\kappa\big(T(U_j)\big)\ \&\ (\forall\sigma\in\Dom(\bar{f}))(
f_\sigma\in {}^{\textstyle{<}\kappa}\kappa)\big\},\]
and $F(U)=\bigcup\limits_{j<\kappa^+} F(U_j)$. Note that $|F(U_j)|\leq
\kappa$ and $|F(U)|\leq\kappa^+$. Choose a function $\psi_1:\kappa^+
\longrightarrow F(U)$ such that  
\[\big(\forall j<\kappa^+\big)\big(\cf(j)=\kappa\ \Rightarrow\ \Rng(\psi_1
\rest j)=F(U_j)\big).\] 
Finally, let $c:\kappa^+\times\kappa^+\longrightarrow\kappa^+$ be a
bijection such that   
\[\big(\forall j<\kappa^+\big)\big(\cf(j)=\kappa\ \Rightarrow\ \Rng(c\rest
(j\times j))=j\big).\]  
Now let $\varphi_{\bar{q}}:\kappa^+\longrightarrow\kappa^+$ be a regressive
function such that for $j<\kappa^+$ of cofinality $\kappa$ we have
\[\varphi_{\bar{q}}(j)=c(\min\{\alpha<\kappa^+: \psi_0(\alpha)=q_j
\rest U_j\},\min\{\alpha<\kappa^+: \psi_1(\alpha)=\bar{f}^{q_j}\rest T(U_j)
\}).\] 
Easily, $\varphi_{\bar{q}}$ is as required.
\end{proof}
\medskip

Now we may complete the proof of Proposition \ref{p633}. Consider the
following strategy $\st$ for Player I in the game $\Gcc(\bQ)$. Suppose that
the players arrived at stage $\alpha>0$ of the play and they have already
constructed a sequence $\langle\bar{q}^\beta,\bar{p}^\beta,\varphi^\beta:
\beta<\alpha\rangle$. Then, for each  $j<\kappa^+$, the sequence $\langle
p^\beta_j:\beta<\alpha\rangle$ is increasing, so Player I can take its least
upper bound $q^\alpha_j$. This determines $\bar{q}^\alpha$ played by Player
I; the function $\varphi^\alpha$ played at this stage is the
$\varphi_{\bar{q}^\alpha}$ given by Claim \ref{cl2}. 

One easily verifies that the strategy $\st$ is a winning one (remember Claim
\ref{cl1}). 
\end{proof}

\begin{theorem}
\label{p3.4}
Suppose $\theta$ and $\iota$ are cardinals such that $\theta=\cf(\theta)<
\iota=\iota^{<\kappa}$. Then the forcing notion $\bQ$ is
$(\theta,\iota,\kappa)$--manageable.  
\end{theorem}

\begin{proof}
For each $\sigma\in\bigcup\limits_{\alpha<\tau}\prod\limits_{\beta<
\alpha}\mu_\beta$ fix a sequence $\eta_\sigma\in\prod\limits_{\alpha<
\tau}\mu_\alpha$ such that $\sigma\subseteq \eta_\sigma$. Let
$\bar{\eta}=\langle \eta_\sigma:\sigma\in\bigcup\limits_{\alpha<\tau}
\prod\limits_{\beta<\alpha}\mu_\beta\rangle$.

Suppose that $N$ is a $(\bQ,\kappa,\iota)$--relevant model such that
$(\bar{\eta},\bar{\mu},\pi)\in N$. 

For a condition $p\in\bQ$ we define conditions $\cl_N^+(p)=q$ and
$\cl_N^-(p)=r$ by  
\begin{itemize}
\item $i^r=i^q=i^p$, 
\item $u^r=(u^p\cap N)\cup\{\eta_\sigma:\sigma\in T(u^p)\cap
N\}$, $u^q=u^p\cup\{\eta_\sigma:\sigma\in T(u^p)\cap N\}$, 
\item $f^r_\sigma=f^p_\sigma$ for $\sigma\in T(u^p)\cap N$,
$f^q_\sigma=f^p_\sigma$ for $\sigma\in T(u^p)$, and $f^q_\sigma(j)=
f^r_\sigma(j)=i^p$ for $\sigma\in T(u^q)\setminus T(u^p)$, $j<i^p$, 
\item $g^r(\eta)=g^p(\eta)$ for $\eta\in u^p\cap N$ and $g^r(\eta)=i^p$ for
$\eta\in u^r\setminus u^p$;\\
$g^q(\eta)=g^p(\eta)$ for $\eta\in u^p$ and $g^q(\eta)=i^p$ for
$\eta\in u^q\setminus u^p$.  
\end{itemize}
Plainly, $\cl^+_N(p), \cl^-_N(p)$ are conditions in $\bQ$ and $\cl^-_N(p)$
belongs to $N$ (remember ${}^{\textstyle{<}\kappa}N\subseteq N$). If $p\in
N$ then also $\cl^+_N(p)\in N$. 

\begin{claim}
\label{cl3}
Suppose that $p\in\bQ$, $q\in N\cap\bQ$ are such that $q\leq p$. Then  
\begin{enumerate}
\item $\cl^+_N(q)=\cl^-_N(q)$, $q\leq\cl_N^-(p)\leq\cl^+_N(p)$, and
$p\leq\cl^+_N(p)$, 
\item if $q'\in N\cap\bQ$ is stronger than $\cl_N^-(p)$, then $q'$ and $p$
are compatible,
\item if $p'\in\bP$ is stronger than $\cl^+_N(p)$, then $\cl^-_N(p)\leq
\cl^-_N(p')$. 
\end{enumerate}
\end{claim}

\begin{proof}[Proof of the Claim]
1)\quad Just check.
\smallskip

\noindent 2)\quad Suppose $\cl^-_N(p)\leq q'\in N\cap\bQ$. Put 
\[i^r=i^{q'},\quad u^r=u^{q'}\cup u^p,\quad g^r(\eta)=\left\{
\begin{array}{ll} g^{q'}(\eta)&\mbox{ if }\eta\in u^{q'},\\
g^p(\eta)&\mbox{ if }\eta\in u^p\setminus u^{q'}\end{array}\right.\quad
\mbox{ and:}\]
if $\sigma\in T(u^{q'})$, then $f^r_\sigma=f^{q'}_\sigma$, and if $\sigma\in
T(u^p)\setminus T(u^{q'})$, then 
\[f^p_\sigma\subseteq f^r_\sigma\quad\mbox{ and }\quad f^r_\sigma(j)=
\sup\{\pi_\eta(j): \sigma\subseteq\eta\in u^p\}+1\quad\mbox{ for }i^p\leq
j<i^{q'}.\] 
Note that if $\eta\in u^p\setminus u^{q'}$, then for some $\alpha<\tau$
we have $\eta\rest\alpha\notin N$ (so $\eta\rest\alpha\notin
T(u^{q'})$). Hence we may easily verify that $r=(i^r,u^r,\bar{f}^r,g^r)\in
\bQ$ and clearly $r$ is stronger than $q'$. To check that it is also
stronger than $p$ it is enough to note that:

if $\eta\in u^p\cap u^{q'}$, then ($\eta\in u^{\cl^-_N(p)}$ and hence)
$g^{q'}(\eta)=g^{\cl^-_N(p)}(\eta)=g^p(\eta)$, and

if $\sigma\in T(u^p)\cap T(u^{q'})$, then ($\sigma\in T(u^{\cl^-_N(p)})$ and
hence) $f^p_\sigma= f^{\cl^-_N(p)}_\sigma\subseteq f^{q'}_\sigma$.
\smallskip

\noindent 3)\quad Note that if $\cl^+_N(p)\leq p'$, then
\[(u^p\cap N)\cup \{\eta_\sigma:\sigma\in T(u^p)\cap N\}\subseteq
u^{p'}\cap N,\] 
so checking the conditions for $\cl^-_N(p)\leq\cl^-_N(p')$ is pretty
straightforward. 
\end{proof}

\begin{claim}
\label{cl4}
Suppose that a sequence $\langle p_\zeta:\zeta<\zeta^*\rangle\subseteq\bQ$
is increasing, $\zeta^*<\kappa$ is a limit ordinal, and $\cl^+_N(p_\zeta)=
p_{\zeta+1}$ for all even $\zeta<\zeta^*$. Let $p^*$ be the least upper
bound of $\langle p_\zeta:\zeta<\zeta^*\rangle$. Then $\cl^-_N(p^*)$ is the
least upper bound of $\langle \cl^-_N(p_\zeta):\zeta<\zeta^*\rangle$. 
\end{claim}

\begin{proof}[Proof of the Claim] It follows from Claim \ref{cl3}(3) that
$\cl^-_N(p_\zeta)\leq \cl^-_N(p^*)$ (for $\zeta<\zeta^*$). To show
that $\cl^-_N(p^*)$ is actually the least upper bound it is enough to
note that
\[i^{\cl^-_N(p^*)}=i^{p^*}=\sup\{i^{p_\zeta}:\zeta<\zeta^*\}=
\sup\{i^{\cl^-_N(p_\zeta)}:\zeta<\zeta^*\},\]
and 
\[\begin{array}{l}
u^{p^*}\cap N=\bigcup\{u^{p_{\zeta+1}}\cap N:\zeta<\zeta^*\ \&\ \zeta\mbox{
even\/}\}=\bigcup\{u^{\cl^-_N(p_\zeta)}:\zeta<\zeta^*\ \&\ \zeta\mbox{
even\/}\},\\  
\{\eta_\sigma\!:\sigma\in T(u^{p^*})\cap N\}=\{\eta_\sigma\!:\sigma\in
T(u^{p_\zeta})\cap N\ \&\ \zeta<\zeta^*\}\subseteq   
\bigcup\{u^{\cl^-_N(p_\zeta)}\!:\zeta<\zeta^*\},
  \end{array}\]
so $u^{\cl^-_N(p^*)}=u^{p^*}\cap N=\bigcup\{u^{\cl^-_N(p_\zeta)}:\zeta<
\zeta^*\}$.  
\end{proof}

Now we may describe a strategy $\st$ for Her in the game $\Game^m(N,\theta,
\bQ)$. Suppose that $i<\theta$ is even and $(q_i,p_i)$ is His move at this
stage of the play (so $q_i\in N\cap \bP$, $q_i\leq p_i\in\bP$). Then $\st$
instructs Her to play $q_{i+1}=\cl^-_N(p_i)$ and $p_{i+1}=\cl^+_N(p_i)$. It
follows from Claim \ref{cl3}(1) that $(q_{i+1},p_{i+1})$ is a legal move. It 
follows from Claims \ref{cl4} and \ref{cl3}(2) that the strategy $\st$ is a
winning one.  

Thus we have shown that $\bQ$ is weakly $(\theta,\iota,\kappa)$--manageable. 
The rest follows from Propositions \ref{p3.3} and \ref{p633}.
\end{proof}

\begin{definition}
\label{d3.5}
We define $\bQ$--names $\dot{f}_\sigma$ (for $\sigma\in\bigcup\limits_{
\alpha<\tau}\prod\limits_{\beta<\alpha}\mu_\beta$) and $\dot{g}$ by 
\[\begin{array}{ll}
\forces_\bQ&\mbox{`` }\dot{f}_\sigma=\bigcup\{f^p_\sigma: p\in\Gamma_\bQ\
\&\ \sigma\in T(u^p)\}\mbox{ '',}\\
\forces_\bQ&\mbox{`` }\dot{g}=\bigcup\{g^p: p\in\Gamma_\bQ\}\mbox{ ''.}
  \end{array}\]
\end{definition}

\begin{proposition}
\label{p3.6}
\begin{enumerate}
\item $\forces_\bQ\mbox{``
}\dot{g}:\prod\limits_{\alpha<\tau}\mu_\alpha\longrightarrow\kappa
\mbox{ ''}$. 
\item For each $\sigma\in\bigcup\limits_{\alpha<\tau}\prod\limits_{
\beta<\alpha}\mu_\beta$ we have\ \  $\forces_\bQ\mbox{`` }\dot{f}_\sigma:
\kappa\longrightarrow\kappa\mbox{ ''}$.
\item For each $\eta\in \prod\limits_{\alpha<\tau}\mu_\alpha$,
\[\forces_\bQ\mbox{`` }\big(\forall j<\kappa\big)\big(\dot{g}(\eta)\leq j\ 
\Rightarrow\ \pi_\eta(j)<\sup\{\dot{f}_{\eta\rest\alpha}(j): \alpha<
\tau\}\big).\]
\end{enumerate}
\end{proposition}

\begin{proof}
For $\eta\in \prod\limits_{\alpha<\tau}\mu_\alpha$ and $i<\kappa$ let 
\[\cI_\eta=\{p\in\bQ:\eta\in u^p\}\quad\mbox{ and }\quad \cI^i=\{p\in\bQ:
i<i^p\}.\] 
We claim that these are open dense subsets of $\bQ$. First, suppose
$\eta\notin u^p$, $p\in\bQ$ and let $i^r=i^p$, $u^r=u^p\cup\{\eta \}$,
$g^p\subseteq g^r$, $g^r(\eta)=i^r$, $f^r_\sigma=f^p_\sigma$ for $\sigma\in
T(u^p)$ and $f^r_{\eta\rest\alpha}(j)=1$ if $\eta\rest\alpha\notin T(u^p)$,
$\alpha<\tau$. Then $r\in\cI_\eta$ is stronger than $p$. (Thus the sets
$\cI_\eta$ are dense.)

Now suppose that $p\in\bP$ is such that $i^p\leq i<\kappa$. Put $i^r=i+1$,
$u^r=u^p$, $g^r=g^p$ and for $\sigma\in T(u^r)$ let $f^r_\sigma\supseteq
f^p_\sigma$ be such that ($\Dom(f^r_\sigma)=i^r$ and) for $j\in i^r\setminus
i^p$ we have $f^r_\sigma(j)=\sup\{\pi_\eta(j):\eta\in u^r\}+1$. This way we
have defined a condition $r\in\bQ$ stronger than $p$ and such that
$r\in\cI^i$. (Thus the sets $\cI^i$ are dense.)

Using the above observation and the definition of the order of $\bQ$ one
easily justifies (1) and (2). (Note also that, as $\bQ$ is $\kappa$--complete,
$\forces_{\bQ}\mbox{`` }\prod\limits_{\alpha<\tau}\mu_\alpha=\big(
\prod\limits_{\alpha<\tau}\mu_\alpha\big)^\bV\mbox{ ''}$.) Then (3)
follows immediately once you note that 
\[p\forces\mbox{`` }g^p\subseteq\dot{g}\ \&\ f^p_\sigma\subseteq
\dot{f}_\sigma\mbox{ ''},\] 
(for $\sigma\in T(u^p)$, $p\in\bQ$); remember Definition \ref{d3.2}(1)(c). 
\end{proof}

\section{The models}

\begin{theorem}
\label{t4.1}
Assume that
\begin{enumerate}
\item[(a)]  $\aleph_0<\kappa=\cf(\kappa)=2^{<\kappa}$,
\item[(b)]  $\mu$ is a cardinal such that $\cf(\mu)<\kappa<\mu<\mu^{\cf(
\mu)}=2^\kappa$, 
\item[(c)] there is an increasing sequence $\bar{\mu}=\langle\mu_\alpha:
\alpha<\cf(\mu)\rangle$ of regular cardinals such that 
\[(\forall\alpha<\cf(\mu))(\kappa\leq\mu_\alpha\leq(\mu_\alpha)^{\cf(\mu)}
\leq\mu_{\alpha+1})
\quad\mbox{ and }\quad \mu=\sup\{\mu_\alpha:\alpha<\cf(\mu)\}.\]
\end{enumerate}
Then there is a forcing notion $\bP$ such that: 
\begin{enumerate}
\item[(i)] $\bP$ has a dense subset of size $2^\kappa$,
\item[(ii)] $\bP$ is $(\theta,\iota,\kappa)$--manageable for all cardinals
$\theta,\iota$ satisfying $\cf(\theta)=\theta<\kappa\leq\iota=\iota^{<
\kappa}$,  
\item[(iii)] $\forces_\bP\mbox{`` }\gd_\kappa^{\cf(\mu)}\leq\mu\mbox{
''}$,
\item[(iv)] if $\cov(\Mkk)>\mu$, then  $\forces_\bP\mbox{`` }\gd_\kappa^{
\cf(\mu)}=\mu<\big(\cov(\Mkk)\big)^\bV\leq\cov(\Mkk)\mbox{ ''}$.
\end{enumerate}
\end{theorem}

\begin{proof}
Assume $\kappa,\mu,\bar{\mu}=\langle\mu_\alpha:\alpha<\cf(\mu)\rangle$ are
as in the assumptions (a)--(c). Note that then also $\prod\limits_{\alpha<
\cf(\mu)}\mu_\alpha=2^\kappa$ (by Tarski's theorem). 
   
The forcing notion $\bP$ is built as the limit of a ${<}\kappa$--support
iteration $\bar{\bQ}=\langle\bP_\xi,\dbQ_\xi:\xi<\kappa^+\rangle$. The names
$\dbQ_\xi$ are defined by induction on $\xi<\kappa^+$ so that 
\begin{enumerate}
\item[$(\alpha)$] $\bP_\xi$ has a dense subset of size $2^\kappa$,
\end{enumerate}
and for all cardinals $\theta,\iota$ satisfying $\cf(\theta)=\theta<\kappa
\leq \iota=\iota^{<\kappa}$,  
\begin{enumerate}
\item[$(\beta)$]  the forcing notion $\bP_\xi$ is
$(\theta,\iota,\kappa)$--manageable and  
\item[$(\gamma)$] $\forces_{\bP_\xi}\mbox{`` $\dbQ_\xi$ is 
$(\theta,\iota,\kappa)$--manageable ''}$.
\end{enumerate}
So suppose that $\bP_\xi$ is already defined (and clauses $(\alpha)$,
$(\beta)$ hold). Then $\bP_\xi$ is $\kappa$--complete $\kappa^+$--cc, and
hence the properties of $\kappa,\mu$ and $\bar{\mu}$ stated in (a)--(c) hold
in $\bV^{\bP_\xi}$. Take a $\bP_\xi$--name 
$\dot{\pi}_\xi$ such that 
\[\forces_{\bP_\xi}\mbox{`` }\dot{\pi}_\xi:\prod\limits_{\alpha<\cf(\mu)}
\mu_\alpha\longrightarrow \bairek\mbox{ is a bijection '',}\]
and let $\dbQ_\xi$ be a $\bP_\xi$--name for the forcing notion $\bQ(
\dot{\pi}_\xi,\bar{\mu},\kappa)$. Then clause $(\gamma)$ holds (remember
Theorem \ref{p3.4}). 

It follows from Proposition \ref{p3.3}(1) that the demand $(\alpha)$ is
preserved at successor stages and it is preserved at limits $\leq\kappa^+$
by the support we use. By Theorem \ref{t1.6}, the clause $(\beta)$ holds for
each $\bP_\xi$ ($\xi\leq\kappa^+$). So our $\bP=\bP_{\kappa^+}$ satisfies
(i)+(ii).  

For $\xi<\kappa^+$ let $\dot{f}^\xi_\sigma$ (for $\sigma\in
\bigcup\limits_{\alpha<\cf(\mu)}\prod\limits_{\beta<\alpha}\mu_\beta$) and 
$\dot{g}^\xi$ be $\bP_{\xi+1}$--names for functions added by $\bQ(
\dot{\pi}_\xi,\bar{\mu},\kappa)$ (see Definition \ref{d3.5}). Then the family 
\[\cF=\{\dot{f}^\xi_\sigma: \xi<\kappa^+\ \&\ \sigma\in
\bigcup\limits_{\alpha<\cf(\mu)}\prod\limits_{\beta<\alpha}\mu_\beta\}\]
is of size $\mu$. Since $\bP$ is $\kappa^+$--cc, for each $\bP$--name
$\dot{h}$ for a member of $\bairek$ there are $\xi<\kappa^+$ and a 
$\bP_\xi$--name $\dot{h}^*$ such that $\forces_{\bP}\mbox{``
}\dot{h}=\dot{h}^*\mbox{ ''}$. Thus using Proposition \ref{p3.6}(3) we get
\[\forces_{\bP}\mbox{`` }(\exists j^*<\kappa)(\exists F\in
\cP_{\cf(\mu)^+}(\cF))(\forall j<\kappa)(j^*\leq j\ \Rightarrow\
\dot{h}(j)<\sup\{f(j):f\in F\})\mbox{ ''.}\]
Now we easily conclude that demand (iii) holds. 

To show (iv) let us assume $\cov(\Mkk)>\mu$. The forcing notion $\bP$ is
$(\aleph_0,\kappa,\kappa)$--manageable, so by Corollary \ref{addedc} and
Proposition \ref{trivial} we have  
\[\forces_\bP\mbox{`` }\mu<\big(\cov(\Mkk)\big)^\bV\leq\cov(\Mkk)\leq
\gd_\kappa\mbox{ ''}.\] 
By (iii) we know $\forces_\bP\mbox{`` }\gd_\kappa^{\cf(\mu)}\leq\mu\mbox{
''}$, but as for each $\alpha<\cf(\mu)$: 
\[\forces_\bP\mbox{`` }(\mu_\alpha)^{\cf(\mu)}\leq\mu_{\alpha+1}<\mu<
\gd_\kappa\mbox{ ''},\]
we immediately get that $\forces_\bP\mbox{`` }\gd_\kappa^{\cf(\mu)}=
\mu\mbox{ ''}$ (remember $\gd_\kappa\leq (\gd_\kappa^{\cf(\mu)})^{\cf(\mu
)}$).
\end{proof}

\begin{corollary}
\label{c4.3}
Assume {\em GCH}. Then there a $\kappa$--complete $\kappa^+$--cc forcing
notion $\bP^*$ such that 
\[\forces_{\bP^*}\mbox{`` }\bgdk=\gd^{\aleph_0}_\kappa=\kappa^{+\omega}\
\mbox{ and }\ \cov(\Mkk)=2^\kappa=\kappa^{+(\omega+1)}\mbox{ ''.}\]
\end{corollary}

\begin{proof}
Let $\bP_0={\mathbb C}_{\kappa^{+(\omega+1)},\kappa}$ be the forcing adding
$\kappa^{+(\omega+1)}$ many Cohen functions in $\bairek$ (with
${<}\kappa$--support). Note that 
\[\forces_{\bP_0}\mbox{`` }
\kappa,\ \mu=\kappa^{+\omega}\mbox{ and }\bar{\mu}=\langle \kappa^{+n}:n<
\omega\rangle\mbox{ are as in Theorem \ref{t4.1}(a)--(c) ''.}\]
Therefore we have a $\bP_0$--name $\dot{\bP}$ for a forcing notion
satisfying \ref{t4.1}(i)--(iv). (Note that $\forces_{\bP_0}\mbox{`` }
\cov(\Mkk)=\kappa^{+(\omega+1)}\mbox{ ''}$, so the assumption of
\ref{t4.1}(iv) holds.)  Let $\bP=\bP_0*\dot{\bP}$.
\end{proof}

Note that Theorem \ref{thone} follows from Corollary \ref{c4.3}, Theorem
\ref{thm1.14}, and the fact that $\cov(\Mkk)\leq\cof(\NSk)$. 

\begin{theorem}
\label{t4.1A}
Assume that (a)--(c) of Theorem \ref{t4.1} hold and 
\begin{enumerate}
\item[(d)] $\nu$ is a regular cardinal such that $\mu<\nu<2^\kappa$.
\end{enumerate}
Then there is a forcing notion $\bP_\nu$ satisfying (i)+(ii) of Theorem
\ref{t4.1} and   
\begin{enumerate}
\item[(iii)$^+_\nu$] $\forces_{\bP_\nu}\mbox{`` }\gd_\kappa^{\tau}=\nu
\mbox{ for every cardinal $\tau$ satisfying $\cf(\mu)\leq\tau<\kappa$ ''}$,  
\item[(iv)$^-$] $\forces_{\bP_\nu}\mbox{`` }\big(\cov(\Mkk)\big)^\bV\leq
\cov(\Mkk)\mbox{ ''}$.
\end{enumerate}
\end{theorem}

\begin{proof}
The forcing notion $\bP_\nu$ is the limit of ${<}\kappa$--support iteration
$\bar{\bQ}=\langle\bP_\xi,\dbQ_\xi:\xi<\nu\rangle$, where $\dbQ_\xi$ are
defined as in the proof of Theorem \ref{t4.1} (so the only difference is the
length of the iteration). As there, $\bP_\nu$ satisfies (i)+(ii) and
$\forces_{\bP_\nu}\mbox{`` }\gd_\kappa^{\cf(\mu)}\leq\nu\mbox{ ''}$.  Since
$\cov(\Mkk)>\kappa$ and $\bP_\nu$ is $(\aleph_0,\kappa,\kappa)$--manageable,
we get (iv)$^-$ (by Corollary \ref{addedc}). To show that (iii)$^+_\nu$
holds, suppose that $\dot{\cF}$ is a $\bP_\nu$--name for a family of
functions in $\bairek$ of size $<\nu$. Then $\dot{\cF}$ is essentially a
$\bP_\xi$--name for some $\xi<\nu$. Since $\bP_{\xi+\kappa}$ adds a subset
of $\kappa$ which is Cohen over $\bV^{\bP_\xi}$, $(\bairek)^{\bV^{\bP_\xi}}$
is not a dominating family in $(\bairek)^{\bV^{\bP_\nu}}$, and hence for any
$\tau<\kappa$  
\[\forces_{\bP_\nu}\mbox{`` $\dot{\cF}$ is not $\tau$--dominating ''}\] 
(remember that each $\bP_\zeta$ is $\kappa$--complete).  
\end{proof}

Now, Theorem \ref{thtwo} follows from Theorem \ref{thm1.14} and the
following Corollary. 

\begin{corollary}
It is consistent, relative to the existence of a cardinal $\nu$ such that
$o(\nu)=\nu^{++}$, that $\bar{\gd}_{\aleph_1}=\aleph_{\omega+1}$ and
$\cov({\bf M}_{\aleph_1,\aleph_1})=\aleph_{\omega+2}$.
\end{corollary}

\begin{proof}
Gitik \cite{Gi89} constructed a model of `` $2^{\aleph_n}<\aleph_\aleph$ for
every $n<\omega$, $\aleph_1^{<\aleph_1}=\aleph_1$ and $2^{\aleph_\omega}=
\aleph_{\omega+2}$ '' from $o(\nu)=\nu^{++}$. Add $\aleph_{\omega+2}$ Cohen
subsets of $\aleph_1$ (with countable support) to the model of Gitik and
then apply Theorem \ref{t4.1A}. 
\end{proof}

\begin{theorem}
\label{t4.5}
Assume that
\begin{enumerate}
\item[(a)] $\kappa=\cf(\kappa)=2^{<\kappa}$, $n<\aleph_0$
\item[(b)] $\mu_0,\mu_1,\ldots,\mu_n$ are cardinals such that 
\[\mu_0>\mu_1>\ldots>\mu_n>\kappa\quad \mbox{ and }\quad \cf(\mu_0)<
\cf(\mu_1)<\ldots<\cf(\mu_n)<\kappa,\] 
\item[(c)] $(\mu_\ell)^{\cf(\mu_\ell)}=2^\kappa$ for $\ell=0,\ldots,n$,
\item[(d)] for $\ell=0,\ldots,n$, there is an increasing sequence
$\bar{\mu}^\ell=\langle\mu^\ell_\alpha:\alpha<\cf(\mu_\ell)\rangle$ of
regular cardinals such that 
\[(\forall\alpha<\cf(\mu_\ell))(\kappa\leq\mu^\ell_\alpha=(
\mu^\ell_\alpha)^{\cf(\mu_\ell)})\quad\mbox{ and }\quad \mu_\ell=\sup\{
\mu^\ell_\alpha:\alpha<\cf(\mu_\ell)\}.\]
\item[(e)] $\cov(\Mkk)=2^\kappa$.
\end{enumerate}
Then there is a forcing notion $\bP$ such that 
\begin{enumerate}
\item[(i)] $\bP$ has a dense subset of size $2^\kappa$,
\item[(ii)] $\bP$ is $(\theta,\iota,\kappa)$--manageable for all cardinals
$\theta,\iota$ satisfying $\cf(\theta)=\theta<\kappa\leq\iota=\iota^{<
\kappa}$,  
\item[(iii)] $\forces_\bP\mbox{`` }\gd_\kappa^{\cf(\mu_\ell)}=\mu_\ell\mbox{
for }\ell=0,\ldots,n\mbox{ and }\cov(\Mkk)=\gd_\kappa=2^\kappa\mbox{ ''}$. 
\end{enumerate}
\end{theorem}

\begin{proof}
Let $A_0,\ldots, A_n$ be a partition of $\kappa^+$ into sets of size
$\kappa^+$. The forcing notion $\bP$ is the limit of a ${<}\kappa$--support
iteration $\bar{\bQ}=\langle\bP_\xi,\dbQ_\xi:\xi<\kappa^+\rangle$ defined
like in the proof of Theorem \ref{t4.1}, but 
\begin{itemize}
\item if $\xi\in A_\ell$, then $\dot{\pi}_\xi$ is a $\bP_\xi$--name for a
bijection from $\prod\limits_{\alpha<\tau_\ell}\mu_\alpha^\ell$ onto
$\bairek$, and $\dbQ_\xi$ is $\bQ(\dot{\pi}_\xi,\bar{\mu}^\ell,\kappa)$.
\end{itemize}
We argue that $\bP$ has the required properties similarly as in Theorem
\ref{t4.1}. 
\end{proof}

\begin{remark}
Of course the assumption (e) in Theorem \ref{t4.5} is not very important: we
may start with adding $2^\kappa$ many Cohen subsets of $\kappa$.
\end{remark}


\def\germ{\frak} \def\scr{\cal} \ifx\documentclass\undefinedcs
  \def\bf{\fam\bffam\tenbf}\def\rm{\fam0\tenrm}\fi 
  \def\defaultdefine#1#2{\expandafter\ifx\csname#1\endcsname\relax
  \expandafter\def\csname#1\endcsname{#2}\fi} \defaultdefine{Bbb}{\bf}
  \defaultdefine{frak}{\bf} \defaultdefine{mathfrak}{\frak}
  \defaultdefine{mathbb}{\bf} \defaultdefine{mathcal}{\cal}
  \defaultdefine{beth}{BETH}\defaultdefine{cal}{\bf} \def\bbfI{{\Bbb I}}
  \def\mbox{\hbox} \def\text{\hbox} \def\om{\omega} \def\Cal#1{{\bf #1}}
  \def\pcf{pcf} \defaultdefine{cf}{cf} \defaultdefine{reals}{{\Bbb R}}
  \defaultdefine{real}{{\Bbb R}} \def\restriction{{|}} \def\club{CLUB}
  \def\w{\omega} \def\exist{\exists} \def\se{{\germ se}} \def\bb{{\bf b}}
  \def\equivalence{\equiv} \let\lt< \let\gt> \def\cite#1{[#1]}
  \def\implies{\Rightarrow}

\end{document}